\documentclass[11pt]{amsart}

\usepackage{amsmath,latexsym,amsbsy,amssymb,enumerate,amsthm}
\usepackage{enumitem}
\usepackage[dvipsnames]{xcolor}
\usepackage{graphicx}
\usepackage{cases}
\usepackage{relsize}
\usepackage{subfigure}
\usepackage{float}
\usepackage{comment}
\usepackage[numbers,sort&compress]{natbib}
\usepackage{tikz}
\usepackage{forest}
\usepackage{tikz-qtree}
\usepackage{amsbsy}
\usepackage[enableskew,stdtext]{youngtab}
\usepackage{ytableau}
\usepackage{curves}
\usepackage{hyperref}
\usepackage[capitalise]{cleveref}
\usepackage{ulem}
\usetikzlibrary{tqft}
\usetikzlibrary{positioning}
\usetikzlibrary{arrows.meta, decorations.pathmorphing,calc,cd}
\def\epsilon{\varepsilon}

\newcommand{\set}[1]{\left\{#1\right\}}
\newcommand{\sett}[2]{\left\{#1 ~\Big| ~ #2\right\}}

\newcommand{\E}[0]{\mathcal{E}}
\newcommand{\Z}[0]{\mathbb{Z}}

\def\zh{\hat{0}}
\def\oneh{\hat{1}}

\newcommand{\precdot}{\prec\mathrel{\mkern-5mu}\mathrel{\cdot}}

\newtheorem{theorem}{Theorem}[section]
\newtheorem{lemma}[theorem]{Lemma}
\newtheorem{corollary}[theorem]{Corollary}
\newtheorem{proposition}[theorem]{Proposition}
\newtheorem{definition}[theorem]{Definition}
\theoremstyle{remark}
\newtheorem{remark}[theorem]{Remark}

\newtheorem{question}[theorem]{Question}

\theoremstyle{definition}

\title{Structure of linear extensions of posets and lexicographic shelling}

\author{Clarice Andorno}
\address{Department of Mathematics, Truman State University, Kirksville, Missouri}
\email{claryando2004@gmail.com}

\author{Stephen Lacina}
\address{Department of Mathematics, Truman State University, Kirksville, Missouri}
\email{slacina@truman.edu}

\begin{document}

\begin{abstract}
We study some structure of linear extensions of a finite poset $P$ by investigating partial orders on the maximal chains of the lattice of order ideals $J(P)$. These partial orders on the maximal chains of $J(P)$ are called maximal chain descent orders (MCDO). MCDOs come from EL-labelings of finite bounded posets and encode the structure of shellings arising from EL-labelings. Every linear extension of $P$ gives an EL-labeling of $J(P)$, and so generates an MCDO. Our main question is: when do different linear extensions of $P$ result in isomorphic MCDOs? We give partial answers to this question. In particular, we prove that linear extensions that are related by an automorphism of $P$ result in isomorphic MCDOs. Additionally, we show that for a certain family of posets, a linear extension $l$ and Sch{\"u}tzenberger's evacuation of $l$ generate isomorphic MCDOs. In this case, the isomorphism comes from a well known isomorphism of the weak order on permutations. %(No AI was used in this work.)
\end{abstract}

\maketitle

\section{Introduction}
There is a deep connection between shelling orders and linear extensions of partially ordered sets (posets). For example, Bj{\"o}rner showed in \cite{coxetercomplexesbjorner1984} that linear extensions of the weak order on a Coxeter group produce shelling orders of the corresponding Coxeter complex. Ardila, Castillo, and Samper showed in \cite{Ardila_Castillo_Samper_2016} that linear extensions of certain partial orders on the bases of a matroid induce shelling orders of related matroid complexes. In \cite{bolognini2023linear}, Bolognini and Sentinelli showed that linear extensions of the Bruhat order (Gale order) on the facets of certain simplicial complexes also give shelling orders of those complexes. They additionally gave a notion of isomorphism of shellings of a simplicial complex based on relabelings of the vertices of the complex. This was built upon by Sentinelli in \cite{sentinelli2025linearly} by introducing a class of simplicial complexes called linearly shellable complexes characterized by admitting shelling orders based on certain linear extensions. 

In the setting of lexicographic shelling of posets, Lacina introduced in \cite{lacina2022maximalchaindescentorders} a partial order on the maximal chains of a finite, bounded poset $P$ with a CL-labeling (or EL-labeling) $\lambda$ called the \textbf{maximal chain descent order (MCDO)}, denoted $\mathbf{MP_{\boldsymbol{\lambda}}}$. (We note that the original notation for an MCDO is $P_{\lambda}(2)$ because of connections with higher Bruhat orders of type A Coxeter groups.) The linear extensions of $MP_{\lambda}$ are precisely the shellings of the order complex of $P$ that have the same restriction map as the lexicographic shellings from $\lambda$. Thus, MCDOs may be thought of as the structure of shellings arising from a CL-labeling (or EL-labeling) of a poset. They also provide a notion of isomorphism of CL-labelings of $P$ which is distinct from Bolognini and Sentinelli's notion of isomorphism of shellings. That is, two CL-labelings of $P$ are isomorphic if their induced MCDOs are isomorphic.

Our present focus is on shellings of finite distributive lattices and their connection to linear extensions. By Birkhoff's Fundamental Theorem of Finite Distributive Lattices in \cite{birkhoffftfdl1937}, every finite distributive lattice is isomorphic to the lattice of order ideals of some finite poset. Let $J(P)$ be such a lattice of ideals of finite poset $P$. Each linear extension of $P$ induces an EL-labeling of $J(P)$ and so lexicographic shellings of the order complex of $J(P)$. These are the $M$-chain EL-labelings of $J(P)$ as a supersolvable lattice given by Stanley in \cite{stanleysupersolvablelats1972}. Thus, each linear extension of $P$ induces an MCDO of $J(P)$,  which by a result in \cite{lacina2022maximalchaindescentorders} is isomorphic to an ideal of the weak order on the symmetric group $S_{|P|}$. This is the structure of linear extensions to which we refer, and this adds to an extensive literature on linear extensions. For instance, Bj{\"o}rner and Wachs studied linear extensions from the perspective of permutation statistics in \cite{permstatslinextsbjornerwachs1991}. Another perspective, first introduced by Sch{\"u}tzenberger in \cite{promotionoriginal} and \cite{evacuationorginal} generalizing his earlier work on Young tableaux and built upon my many others (see Stanley \cite{stanley2009} for a survey), is to define operators on linear extensions called promotion and evacuation and study the algebraic structure of these operators. Promotion and evacuation are not easy to understand for general posets, but in certain cases, we provide a connection with MCDOs through partial answers to our motivating question.

\begin{question}\label{quest:isom mcdos?}
When do distinct linear extensions of $P$ produce isomorphic MCDOs of $J(P)$, that is, when do distinct linear extensions of $P$ produce isomorphic EL-labelings of $J(P)$?
\end{question}

We give partial answers to this question beginning with our first main result which we prove as \cref{thm:technical poset automorphisms give isom mcdos}. 

\begin{theorem}\label{thm:poset automorphisms give isom mcdos}
    Let $P$ be a finite poset with linear extension $l$. Let $f:P\to P$ be a poset automorphism. Then the linear extensions $l$ and $l\circ f$ induce isomorphic MCDOs of $J(P)$.
\end{theorem}

\cref{thm:poset automorphisms give isom mcdos} suggests thinking of these isomorphisms of MCDOs of $J(P)$ as a sort of symmetry of $P$. Although we must think of this as a ``generalized" symmetry of $P$ because different linear extensions of $P$ may give isomorphic MCDOs even if the linear extensions are not related by an automorphism of $P$. See \cref{fig:mcdo} and \cref{fig;isolinext} for such examples. These MCDO isomorphisms are actions an the maximal chains of $J(P)$ in contrast to work such as \cite{actiononideals}, \cite{chainhomomesy}, and \cite{minusculeperiodicity} which study actions on ideals of a poset. 

In our second main theorem we expound on such MCDO isomorphisms and illustrate a connection to evacuation by proving that if $P$ is a disjoint union of chains, then a linear extension and its evacuation induce isomorphic MCDOs.

\begin{theorem}\label{thm:evacuation of disjoint union of chains}
    Let $P$ be a disjoint union of finite chains. Let $l$ be a linear extension with $\varepsilon l$ the evacuation of $l$. Then $l$ and $\varepsilon l$ induce isomorphic MCDOs of $J(P)$.
\end{theorem}

We prove this as \cref{thm:evacuation and disjoint chains}. There are two central components of the  proof the first of which is \cref{thm:evacuation formula for disjoint chains} which gives a formula for the evacuation in the case of a disjoint union of chains. This formula was previously proven independently by Nguyen in \cite{cactusgroup}, but we include a proof for completeness and for convenience since their results are expressed in somewhat more generality than we require. The second component is observing that this formula coincides with the restriction of a well-known isomorphism of weak order on $S_n$.

Lastly, our third theorem, which we prove as \cref{thm:stacking preserves mcdo isom}, shows that ordinal sums (stacking posets) preserve isomorphism of MCDOs. For this, we must combine linear extensions in a ``natural" way given in \cref{rmk:properties of ordinal sums}.

\begin{theorem}\label{thm:stacking preserves isom}
    Let $P$ and $Q$ be finite posets. Let $k$ and $k'$ be linear extensions of $P$ that give rise to isomorphic MCDOs of $J(P)$. Similarly, let $l$ and $l'$ be linear extensions of $Q$ producing isomorphic MCDOs of $J(Q)$. Then the natural linear extensions $k\oplus l$ and $k'\oplus l'$ of the ordinal sum $P\oplus Q$ produce isomorphic MCDOs of $J(P\oplus Q)$.
\end{theorem}

This can be easily extended to show that ordinal sums preserve the property of evacuation inducing isomorphisms of MCDOs. Combined with \cref{thm:evacuation and disjoint chains} this might lead one to wonder if evacuation always induces MCDO isomorphisms for series parallel posets. However, the example in \cref{fig:series parallel example} shows that this is not the case. This example shows that disjoint union of posets does not preserve isomorphisms induced by evacuation.

The remainder of this paper has the following structure: \cref{sec:background} provides some necessary background on posets, evacuation, lexicographic shelling, and MCDOs. Then \cref{sec:results} contains our main theorems and their proofs.

\section{Background}\label{sec:background}
\subsection{Posets}
For general poset terminology, notation, and background we refer to Stanley \cite{ec1}. Let $P$ be a finite poset with relation $\leq$. We denote cover relations in $P$ as $x\lessdot y$ when $x$ is covered by $y$. For $x\leq y$ in $P$, the \textbf{closed interval between $x$ and $y$} is the set $[x,y] =\left\{z\in P: x\leq z\leq y \right\}$ which inherits a partial order from $P$. We say that $P$ is \textbf{bounded} if it has a unique maximal element, denoted $\oneh$, and a unique minimal element, denoted $\zh$. A \textbf{chain} of $P$ is a totally ordered subset of $P$. A \textbf{maximal chain} is a chain that is not properly contained in any other chain. The \textbf{length} of a chain with $k+1$ elements is $k$. $P$ is said to be \textbf{graded} if every maximal chain of $P$ has the same length. We use $[n]$ to denote the chain of positive integers at most $n$ with the usual total order. A \textbf{linear extension} of $P$ is a bijective function $l:P\to \left[|P|\right]$ such that $l(x)<l(y)$ whenever $x< y$ in $P$. Posets $P$ and $Q$ are said to be \textbf{isomorphic} if there is a function $f:P\to Q$ such that $f$ is a bijection and $x\leq y$ in $P$ if and only if $f(x)\leq f(y)$ in $Q$. Such a function $f$ is called a \textbf{poset isomorphism}. For posets $P$ and $Q$, the \textbf{disjoint union} $P+Q$ is the partial order on the disjoint union of $P$ and $Q$ given by $x\leq y$ if and only if $x\leq y$ in $P$ or $x\leq y$ in $Q$. Intuitively, $P+Q$ is $P$ next to $Q$. In contrast, the \textbf{ordinal sum} $P\oplus Q$ is the partial order on the disjoint union of $P$ and $Q$ with $x\leq y$ if and only if $x\leq y$ in $P$, or $x\leq y$ in $Q$, or $x\in P$ and $y\in Q$. In spirit, $P\oplus Q$ is $Q$ stacked on top of $P$.

An \textbf{order ideal (ideal)} of $P$ is a subset $I\subseteq P$ such that if $y\in I$ and $x\leq y$ in $P$, then $x\in P$, that is a downwardly closed subset of $P$. The set $J(P)$ of ideals of $P$ ordered by containment is called the \textbf{lattice of ideals of $P$}. $J(P)$ has $\zh=\emptyset$ and $\oneh=P$. By Birkhoff's Fundamental Theorem of Finite Distributive Lattices in \cite{birkhoffftfdl1937}, every finite distributive lattice is isomorphic to the lattice of ideals of some finite poset. See \cref{fig;poset stuff examples} for an example of a poset, a linear extension, and the lattice of ideals.

\subsection{EL-labelings and Maximal Chain Descent Orders}

Each poset $P$ has an associated (abstract) simplicial complex $\Delta(P)$ called the \textbf{order complex of $P$} whose $i$-dimensional faces are the length $i$ chains of $P$. The facets (maximal faces by containment) of $\Delta(P)$ are precisely the maximal chains of $P$. We say a poset has a topological property if its order complex has that property. See Wachs \cite{posettoptools} for a thorough survey of poset topology. A \textbf{shelling} of a simplicial complex $\Delta$ is a total order on the facets $F_1,F_2,\dots, F_t$ of $\Delta$ such that for each $1<j\leq t$, $\overline{F_j} \cap \left(\bigcup_{1\leq i\leq j}\overline{F_i} \right)$ is pure (all facets have the same dimension) and is a codimension one subcomplex of $\overline{F_j}$. Here $\overline{F_i}$ is the subcomplex generated by $F_i$. If $\Delta$ admits a shelling, then $\Delta$ is said to be \textbf{shellable}. We say $P$ is shellable if $\Delta(P)$ is shellable. A shellable simplicial complex is homotopy equivalent to a wedge of top dimensional spheres, is Cohen-Macaulay, and posesses several other nice combinatorial, topological, and algebraic properties.

Lexicographic shelling is a family of techniques to give shellings of order complexes of finite posets based on assigning labels to the cover relations of the poset and then ordering maximal chains based on lexicographic (dictionary) order of associated label sequences. Lexicographic shelling was first introduced by Bj{\"o}rner in \cite{shellblposets} and subsequently expanded by Bj{\"o}rner and Wachs in \cite{bruhatshellbjornerwachs1982}, \cite{lexshellposetsbjornerwachs1983}, \cite{shllblnonpureIbjornerwachs1996}, and \cite{shllblnonpureIIbjornerwachs1997} along with others. More formally, an \textbf{edge labeling of $P$} is a function $\lambda:\E(P) \to \Lambda$ from the set $\E(P)$ of cover relations of $P$ to a poset $\Lambda$. In this paper, $\Lambda=\Z$ with its usual total order. For a maximal chain $m:x_0\lessdot x_1\lessdot\dots \lessdot x_n$, the \textbf{label sequence of $m$} denoted $\lambda(m)$ is the sequence $\lambda(x_0,x_1),\lambda(x_1,x_2),\dots, \lambda(x_{n-1},x_n)$. We say that $m$ is \textbf{ascending} if $\lambda(m)$ is non-decreasing with respect to $\Lambda$. We say that $x\lessdot y\lessdot z$ is an \textbf{ascent} if $\lambda(x,y)\leq \lambda(y,z)$ and a \textbf{descent} if $\lambda(x,y)>\lambda(y,z)$ with respect to $\Lambda$. We say that maximal chain $m$ \textbf{lexicographically precedes} maximal chain $m'$ if $\lambda(m)$ lexicographically precedes $\lambda(m')$, that is if in the first index where $\lambda(m)$ and $\lambda(m')$ differ, $\lambda(m)$ has a smaller label than $\lambda(m')$ with respect to $\Lambda$. 

\begin{definition}[Section 2 \cite{shellblposets}]
    An edge labeling $\lambda$ of a finite, bounded poset $P$ is called an \textbf{EL-labeling} if for each $x<y$, there is a unique ascending maximal chain $m$ in $[x,y]$ with respect to $\lambda$ and $m$ lexicographically precedes every other maximal chain of $[x,y]$.
\end{definition}

Theorem 2.3 in \cite{shellblposets} says that if finite bounded poset $P$ admits an EL-labeling $\lambda$, then any total order on the maximal chains that is compatible with lexicographic order on the label sequences induces a shelling order of $\Delta(P)$. Thus, $P$ is shellable. 

\begin{figure}[h]
\begin{center}

\begin{tikzpicture}[very thick,scale=0.98]
        \node at (0,0) {};
        \draw [fill] (0,0) circle [radius=0.05];
        \node [below] at (0,0) {$a$};
        \node at (-1,1) {};
        \draw [fill] (-1,1) circle [radius=0.05];
        \node [left] at (-1,1) {$b$};
        \node at (1,1) {};
        \draw [fill] (1,1) circle [radius=0.05];
        \node [right] at (1,1) {$c$};
        \node at (-1,2) {};
        \draw [fill] (-1,2) circle [radius=0.05];
        \node [left] at (-1,2) {$d$};
        \node at (1,2) {};
        \draw [fill] (1,2) circle [radius=0.05];
        \node [right] at (1,2) {$e$};
        \node at (1,3) {};
        \draw [fill] (1,3) circle [radius=0.05];
        \node [right] at (1,3) {$f$};
        \draw (0,0) -- (-1,1) -- (-1,2);
        \draw (0,0) -- (1,1) -- (1,2) -- (1,3);
    \end{tikzpicture}
     \hspace{0.5cm}
 \begin{tikzpicture}[very thick,scale=0.98]
        \node at (0,0) {};
        \draw [fill] (0,0) circle [radius=0.05];
        \node [below] at (0,0) {$1$};
        \node at (-1,1) {};
        \draw [fill] (-1,1) circle [radius=0.05];
        \node [left] at (-1,1) {$2$};
        \node at (1,1) {};
        \draw [fill] (1,1) circle [radius=0.05];
        \node [right] at (1,1) {$3$};
        \node at (-1,2) {};
        \draw [fill] (-1,2) circle [radius=0.05];
        \node [left] at (-1,2) {$4$};
        \node at (1,2) {};
        \draw [fill] (1,2) circle [radius=0.05];
        \node [right] at (1,2) {$5$};
        \node at (1,3) {};
        \draw [fill] (1,3) circle [radius=0.05];
        \node [right] at (1,3) {$6$};
        \draw (0,0) -- (-1,1) -- (-1,2);
        \draw (0,0) -- (1,1) -- (1,2) -- (1,3);
    \end{tikzpicture}
    \hspace{0.1cm}
    \begin{tikzpicture}[very thick,scale=1.2]
        \node at (0,0) {$\emptyset$};
        \node at (0,1) {$\{a\}$};
        \node at (-1,2) {$\{a,b\}$};
        \node at (1,2) {$\{a,c\}$};
        \node at (-2,3) {$\{a,b,d\}$};
        \node at (0,3) {$\{a,b,c\}$};
        \node at (2,3) {$\{a,c,e\}$};
        \node at (-1,4) {$\{a,b,c,d\}$};
        \node at (1,4) {$\{a,b,c,e\}$};
        \node at (3,4) {$\{a,c,e,f\}$};
        \node at (0,5) {$\{a,b,c,d,e\}$};
        \node at (2,5) {$\{a,b,c,e,f\}$};
        \node at (1,6) {$\{a,b,c,d,e,f\}$};
        \draw (0,0.2) --node[left,blue] {\scriptsize $1$} (0,0.8);
        \draw (0,1.2) --node[left,blue] {\scriptsize $2$} (-1,1.8);
        \draw (0,1.2) --node[left,blue] {\scriptsize$3$} (1,1.8);
        \draw (-1,2.2) --node[left,blue] {\scriptsize$4$} (-2,2.8);
        \draw (-1,2.2) --node[left,blue] {\scriptsize$3$} (0,2.8);
        \draw (1,2.2) --node[left,blue] {\scriptsize$2$} (0,2.8);
        \draw (1,2.2) --node[left,blue] {\scriptsize$5$} (2,2.8);
        \draw (-2,3.2) --node[left,blue] {\scriptsize$3$} (-1,3.8);
        \draw (0,3.2) --node[left,blue] {\scriptsize$5$} (1,3.8);
        \draw (0,3.2) --node[left,blue] {\scriptsize$4$} (-1,3.8);
        \draw (2,3.2) --node[left,blue] {\scriptsize$2$} (1,3.8);
        \draw (2,3.2) --node[left,blue] {\scriptsize$6$} (3,3.8);
        \draw (-1,4.2) --node[left,blue] {\scriptsize$5$} (0,4.8);
        \draw (1,4.2) --node[left,blue] {\scriptsize$4$} (0,4.8);
        \draw (1,4.2) --node[left,blue] {\scriptsize$6$} (2,4.8);
        \draw (3,4.2) --node[left,blue] {\scriptsize$2$} (2,4.8);
        \draw (0,5.2) --node[left,blue] {\scriptsize$6$} (1,5.8);
        \draw (2,5.2) --node[left,blue] {\scriptsize$4$} (1,5.8);
    \end{tikzpicture}
    \caption{A poset (left) with a linear extension (middle) and the resulting EL-labeling of its lattice of ideals (right).}\label{fig;poset stuff examples}
\end{center}
\end{figure}
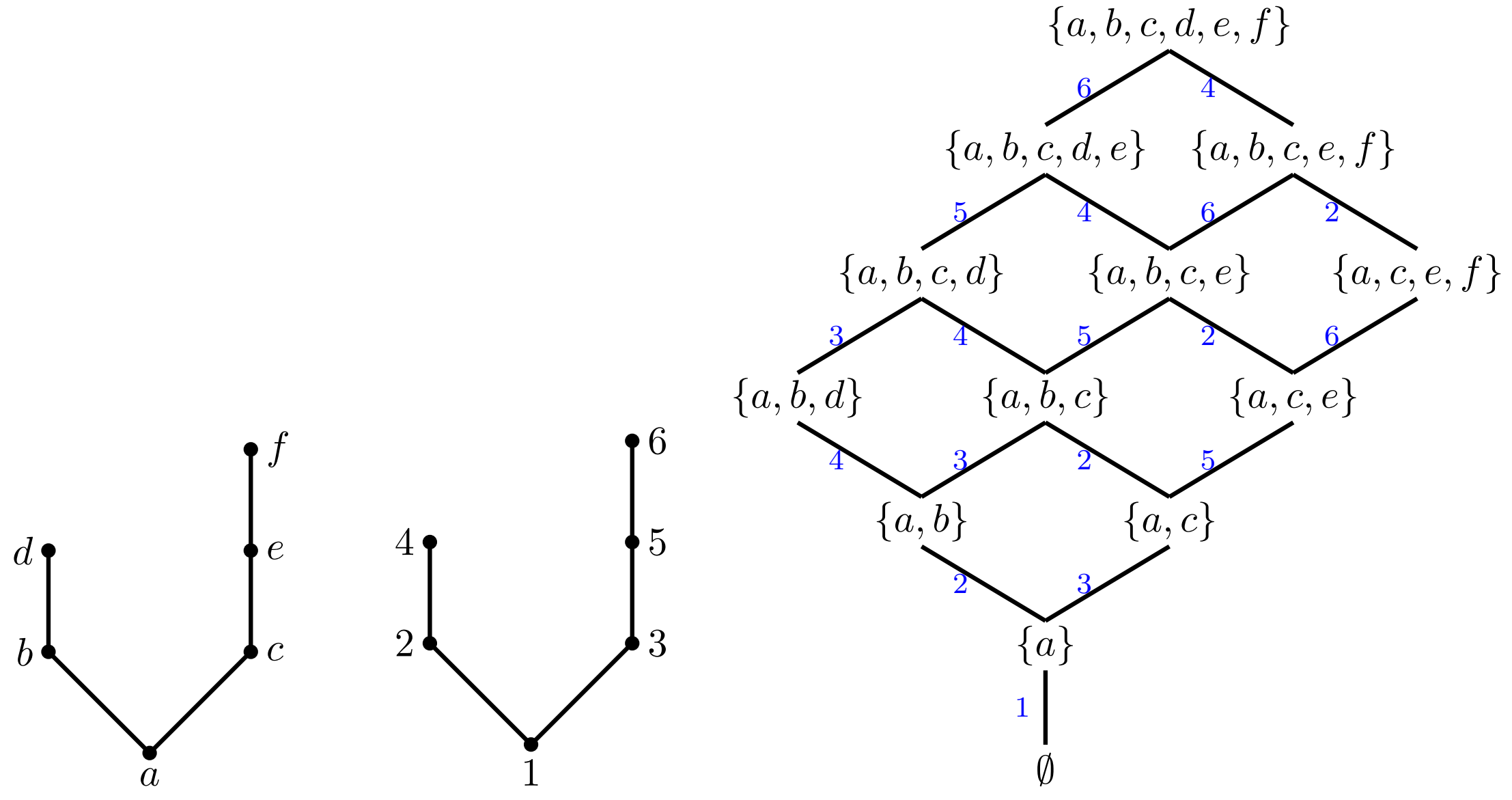

For a finite poset $P$, each linear extension $l$ of $P$ induces an EL-labeling $\lambda_l$ of $J(P)$ as follows. For any two ideals $I_1$ and $I_2$ of $P$ with $I_1 \lessdot I_2$ in $J(P)$, we have $I_2\setminus I_1=\{x\}$ for some $x$ in $P$. We label the cover relation by $\lambda_l (I_1 \lessdot I_2) = l(x)$. \cref{fig;poset stuff examples} gives an example of this type of labeling. These labelings essentially appear in \cite{shellblposets} though these labelings were known earlier outside the context of lexicographic shelling, and $J(P)$ was known to be shellable earlier by other means. These EL-labelings are precisely the $M$-chain EL-labelings of $J(P)$ as a supersolvable lattice due to Stanley in \cite{stanleysupersolvablelats1972}. (They were actually a motivation for the general $M$-chain EL-labelings). 

\begin{figure}[h]
\begin{center}
    
    \begin{tikzpicture}[very thick,scale=1]
        \node at (0,0) { $\emptyset$};
        \node at (0,1) { $\{a\}$};
        \node at (-1,2) { $\{a,b\}$};
        \node at (1,2) { $\{a,c\}$};
        \node at (-2,3) { $\{a,b,d\}$};
        \node at (0,3) {$\{a,b,c\}$};
        \node at (2,3) { $\{a,c,e\}$};
        \node at (-1,4) { $\{a,b,c,d\}$};
        \node at (1,4) { $\{a,b,c,e\}$};
        \node at (3,4) {$\{a,c,e,f\}$};
        \node at (0,5) { $\{a,b,c,d,e\}$};
        \node at (2,5) { $\{a,b,c,e,f\}$};
        \node at (1,6) { $\{a,b,c,d,e,f\}$};
        \draw[orange] (0,0.2) --node[left,blue] {\scriptsize $1$} (0,0.8);
        \draw[purple] (0.1,0.2) -- (0.1,0.8);
        \draw[orange] (0,1.2) --node[left,blue] {\scriptsize$2$} (-1,1.8);
        \draw[purple] (0,1.2) --node[left,blue] {\scriptsize$3$} (1,1.8);
        \draw (-1,2.2) --node[left,blue] {\scriptsize$4$} (-2,2.8);
        \draw [orange] (-1,2.2) --node[left,blue] {\scriptsize$3$} (0,2.8);
        \draw [purple](1,2.2) --node[left,blue] {\scriptsize$2$} (0,2.8);
        \draw (1,2.2) --node[left,blue] {\scriptsize$5$} (2,2.8);
        \draw (-2,3.2) --node[left,blue] {\scriptsize$3$} (-1,3.8);
        \draw[purple] (0.1,3.1) -- (1.1,3.7);
        \draw[orange] (0,3.2) --node[left,blue] {\scriptsize$5$} (1,3.8);
        \draw (0,3.2) --node[left,blue] {\scriptsize$4$} (-1,3.8);
        \draw (2,3.2) --node[left,blue] {\scriptsize$2$} (1,3.8);
        \draw (2,3.2) --node[left,blue] {\scriptsize$6$} (3,3.8);
        \draw (-1,4.2) --node[left,blue] {\scriptsize$5$} (0,4.8);
        \draw[purple] (1.1,4.3) -- (0.1,4.9);
        \draw[orange] (1,4.2) --node[left,blue] {\scriptsize$4$} (0,4.8);
        \draw (1,4.2) --node[left,blue] {\scriptsize$6$} (2,4.8);
        \draw (3,4.2) --node[left,blue] {\scriptsize$2$} (2,4.8);
        \draw[purple] (0.1,5.1) -- (1.1,5.7);
        \draw[orange] (0,5.2) --node[left,blue] {\scriptsize$6$} (1,5.8);
        \draw (2,5.2) --node[left,blue] {\scriptsize$4$} (1,5.8);
    \end{tikzpicture}
    \hspace{0.1cm}
    \begin{tikzpicture}[very thick,scale=1]
        \node[blue] at (0,0) {\scriptsize$123456$};
        \node[blue] at (0,1) {\scriptsize$132456$};
        \node[blue] at (-2,1) {\scriptsize$124356$};
        \node[orange] at (2,1) {\scriptsize$123546$};
        
        \node[purple] at (0,2) {\scriptsize$132546$};

        \node[blue] at (0,3) {\scriptsize$135246$};
        \node[blue] at (4,3) {\scriptsize$132564$};
        \node[blue] at (4,2) {\scriptsize$123564$};

        \node[blue] at (2,4) {\scriptsize$135264$};

        \node[blue] at (2,5) {\scriptsize$135624$};

        \draw (0,0.2) -- (-2,0.8);
        \draw (0,0.2) -- (0,0.8);
        \draw (0,1.2) -- (0,1.8);
        \draw (0,2.2) -- (0,2.8);
        \draw (0,3.2) -- (2,3.8);
        \draw (2,4.2) -- (2,4.8);
        \draw (0,0.2) -- (2,0.8);
        \draw (2,1.2) -- (4,1.8);
        \draw (4,2.2) -- (4,2.8);
        \draw (4,3.2) -- (2,3.8);
        \draw (2,1.2) -- (0,1.8);
        \draw (0,2.2) -- (4,2.8);

    \end{tikzpicture}
    
\end{center}
\caption{An EL-labeling of $J(P)$ from \cref{fig;poset stuff examples} (left) and its MCDO (right).}\label{fig:mcdo}
\end{figure}
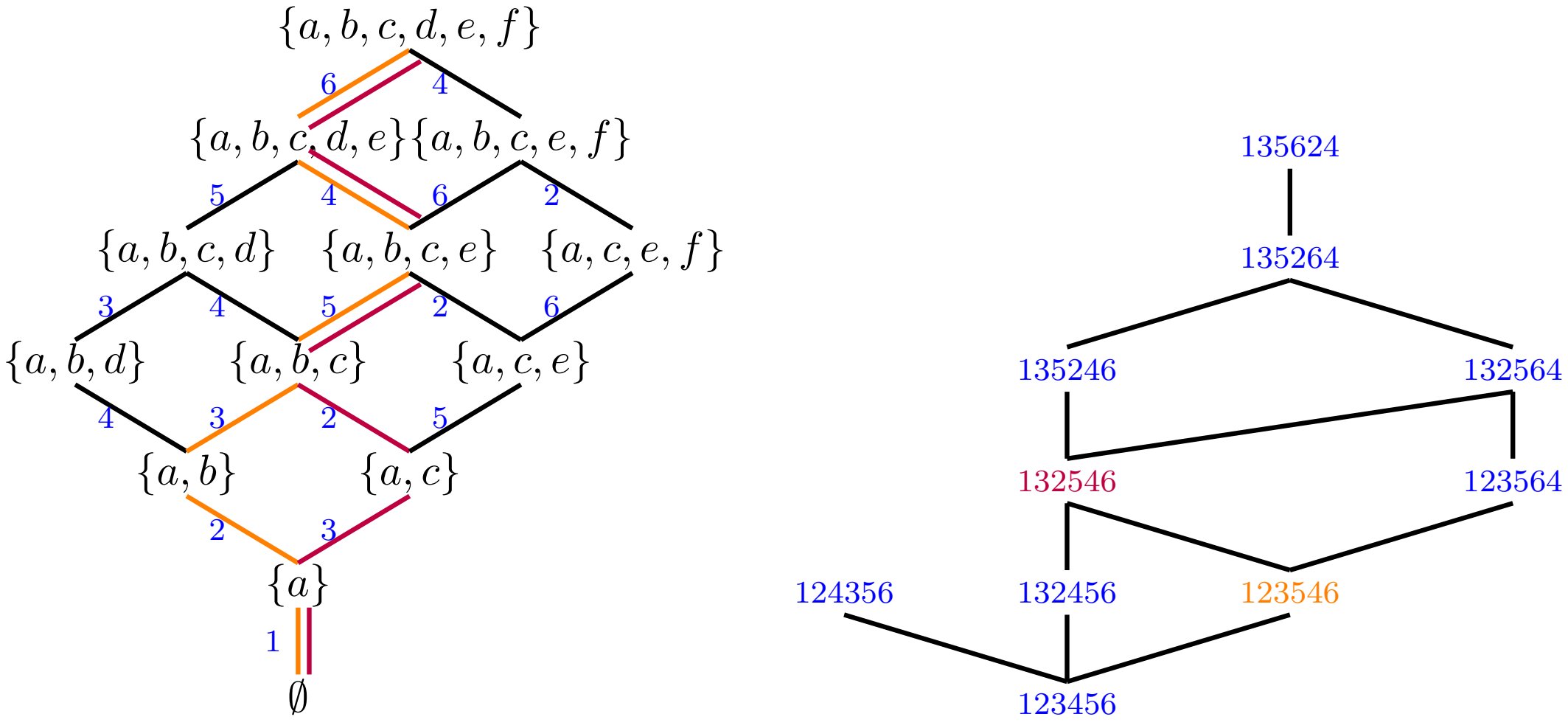

In \cite{lacina2022maximalchaindescentorders}, Lacina introduced a partial order on the maximal chains of a finite, bounded poset $P$ with an EL-labeling $\lambda$ (and more generally a CL-labeling). This partial order on maximal chains encodes the structure of the shellings induced by the EL-labeling. These partial orders are the central object of this paper and give a notion of isomorphism of EL-labelings. We present the definition here for EL-labelings of graded posets since it is simpler and that is all we need, though the original construction does not require gradedness. Say maximal chains $m:x_0\lessdot x_1\lessdot\dots \lessdot x_n$ and $m':x'_0\lessdot x'_1\lessdot\dots \lessdot x'_n$ \textbf{differ in one position} if there is some $1\leq i \leq n-1$ such that $x_i\neq x'_i$ while $x_j=x'_j$ for all $j\neq i$. Further, we say that there is a \textbf{diamond move from $m$ to $m'$} if $m$ and $m'$ differ in one position $i$ and in the interval $\left[x_{i-1},x_{i+1}\right]$, we have that $m$ is the unique ascending chain with respect to $\lambda$. For example, in \cref{fig:mcdo}, the maximal chains $M:\emptyset\lessdot \{a\}\lessdot \{a,b\}\lessdot \{a,b,c\}\lessdot \{a,b,c,e\}\lessdot \{a,b,c,d,e\}\lessdot \{a,b,c,d,e,f\}$ and $M':\emptyset\lessdot \{a\}\lessdot \{a,c\}\lessdot \{a,b,c\}\lessdot \{a,b,c,e\}\lessdot \{a,b,c,d,e\}\lessdot \{a,b,c,d,e,f\}$ differ in one position and there is a diamond move from $M$ to $M'$.

\begin{definition}[Definition 3.6 \cite{lacina2022maximalchaindescentorders}]\label{def:mcdo def}
   Let $P$ be a finite, bounded poset admitting an EL-labeling $\lambda$. The \textbf{maximal chain descent order of $P$ induced by $\lambda$ (MCDO)} is the partial order on the maximal chains of $P$ given by the reflexive and transitive closure of all diamond moves. We denote the order relation on the MCDO induced by $\lambda$ as $\preceq_{\lambda}$ with cover relations denoted $\precdot_{\lambda}$.
\end{definition}

\begin{remark}\label{rmk:mcdo properties}
    We remark here about notation and cover relations in MCDOs. Given the setup of \cref{def:mcdo def}, the original notation for the MCDO of $P$ induced by $\lambda$ was $P_{\lambda}(2)$ because of a connection with higher Bruhat orders of type A Coxeter groups. Presently we are focused solely on the MCDOs of $J(P)$ induced by $\lambda_l$ for some linear extension $l$ of the finite poset $P$. So, for the sake of brevity, we denote such an MCDO by $MP_{\lambda_l}$. 

    In general, the cover relations in an MCDO are somewhat subtle in the sense that not all diamond moves necessarily give cover relations. However, for MCDOs of the form $MP_{\lambda_l}$ the cover relations are exactly the diamond moves by Corollary 5.2 in \cite{lacina2022maximalchaindescentorders}. 
\end{remark}

For the EL-labelings $\lambda_l$, it is clear that the maximal chains of $J(P)$ have distinct label sequences. Thus, when we draw the Hasse diagram of an MCDO $MP_{\lambda_l}$, we may represent a maximal chain by its label sequence. \cref{fig:mcdo} gives an example of an MCDO from a linear extension EL-labeling of a lattice of ideals.

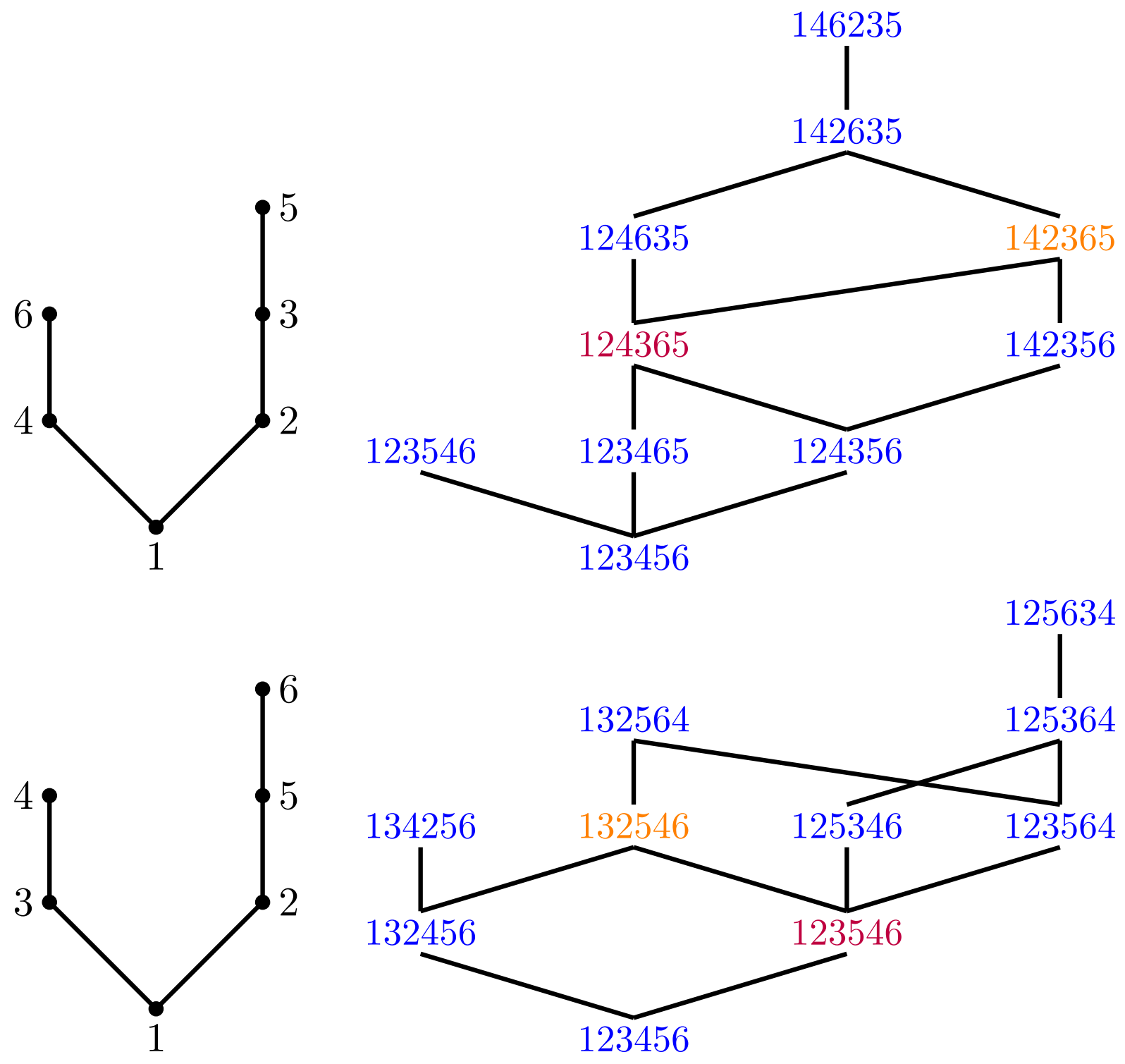
\begin{figure}[H]
 \begin{center}
     \begin{tikzpicture}[very thick,scale=1]
        \node at (0,0) {};
        \draw [fill] (0,0) circle [radius=0.05];
        \node [below] at (0,0) {$1$};
        \node at (-1,1) {};
        \draw [fill] (-1,1) circle [radius=0.05];
        \node [left] at (-1,1) {$4$};
        \node at (1,1) {};
        \draw [fill] (1,1) circle [radius=0.05];
        \node [right] at (1,1) {$2$};
        \node at (-1,2) {};
        \draw [fill] (-1,2) circle [radius=0.05];
        \node [left] at (-1,2) {$6$};
        \node at (1,2) {};
        \draw [fill] (1,2) circle [radius=0.05];
        \node [right] at (1,2) {$3$};
        \node at (1,3) {};
        \draw [fill] (1,3) circle [radius=0.05];
        \node [right] at (1,3) {$5$};
        \draw (0,0) -- (-1,1) -- (-1,2);
        \draw (0,0) -- (1,1) -- (1,2) -- (1,3);
    \end{tikzpicture}
    \hspace{0.1cm}
    \begin{tikzpicture}[very thick,scale=1]
        \node[blue] at (0,0) {\small$123456$};
        \node[blue] at (2,1) {\small$124356$};
        \node[blue] at (-2,1) {\small$123546$};
        \node[blue] at (0,1) {\small$123465$};
        
        \node[purple] at (0,2) {\small$124365$};

        \node[blue] at (0,3) {\small$124635$};
        \node[orange] at (4,3) {\small$142365$};
        \node[blue] at (4,2) {\small$142356$};

        \node[blue] at (2,4) {\small$142635$};

        \node[blue] at (2,5) {\small$146235$};

        \draw (0,0.2) -- (-2,0.8);
        \draw (0,0.2) -- (0,0.8);
        \draw (0,1.2) -- (0,1.8);
        \draw (0,2.2) -- (0,2.8);
        \draw (0,3.2) -- (2,3.8);
        \draw (2,4.2) -- (2,4.8);
        \draw (0,0.2) -- (2,0.8);
        \draw (2,1.2) -- (4,1.8);
        \draw (4,2.2) -- (4,2.8);
        \draw (4,3.2) -- (2,3.8);
        \draw (2,1.2) -- (0,1.8);
        \draw (0,2.2) -- (4,2.8);

    \end{tikzpicture}
    \hspace{10cm}
   \begin{tikzpicture}[very thick,scale=1]
        \node at (0,0) {};
        \draw [fill] (0,0) circle [radius=0.05];
        \node [below] at (0,0) {$1$};
        \node at (-1,1) {};
        \draw [fill] (-1,1) circle [radius=0.05];
        \node [left] at (-1,1) {$3$};
        \node at (1,1) {};
        \draw [fill] (1,1) circle [radius=0.05];
        \node [right] at (1,1) {$2$};
        \node at (-1,2) {};
        \draw [fill] (-1,2) circle [radius=0.05];
        \node [left] at (-1,2) {$4$};
        \node at (1,2) {};
        \draw [fill] (1,2) circle [radius=0.05];
        \node [right] at (1,2) {$5$};
        \node at (1,3) {};
        \draw [fill] (1,3) circle [radius=0.05];
        \node [right] at (1,3) {$6$};
        \draw (0,0) -- (-1,1) -- (-1,2);
        \draw (0,0) -- (1,1) -- (1,2) -- (1,3);
    \end{tikzpicture}
    \hspace{0.1cm}
     \begin{tikzpicture}[very thick,scale=1]
        \node[blue] at (0,0) {\small$123456$};
        %\node[blue] at (0,1) {$132456$};
        \node[blue] at (-2,1) {\small$132456$};
        \node[blue] at (-2,2) {\small$134256$};
        \node[purple] at (2,1) {\small$123546$};
        \node[blue] at (2,2) {\small$125346$};
        
        \node[orange] at (0,2) {\small$132546$};

        \node[blue] at (0,3) {\small$132564$};
        \node[blue] at (4,3) {\small$125364$};
        \node[blue] at (4,2) {\small$123564$};
        \node[blue] at (4,4) {\small$125634$};

        \draw (0,0.2) -- (-2,0.8);
        \draw (-2,1.2) -- (-2,1.8);
        \draw (-2,1.2)-- (0,1.8);
        \draw (4,2.2) -- (0,2.8);
        \draw (0,2.2) -- (0,2.8);
        \draw (0,0.2) -- (2,0.8);
        \draw (2,1.2) -- (4,1.8);
        \draw (4,2.2) -- (4,2.8);
        \draw (2,1.2) -- (0,1.8);
        \draw (4,3.2) -- (4,3.8);
        \draw (2,1.2) -- (2,1.8);
        \draw (2,2.2) -- (4,2.8);
    \end{tikzpicture}
    \caption{Two more linear extensions of $P$ from \cref{fig;poset stuff examples} and their corresponding MCDOs}\label{fig;isolinext}
 \end{center}
\end{figure}

Following that in \cref{fig;isolinext} are two more examples of MCDOs resulting from different linear extensions of the same poset as in \cref{fig:mcdo}. The top is isomorphic to that in \cref{fig:mcdo} while the bottom is not. Furthermore, the linear extension used for the top is the evacuation of the linear extension used in \cref{fig:mcdo}. Thus, these three examples motivate our main questions.

\subsection{Evacuation}
Sch{\"u}tzenberger introduced promotion and evacuation, two operators on linear extensions of a poset, in \cite{promotionoriginal} and \cite{evacuationorginal}. He originally defined these for Young Tableaux and then generalized them to linear extensions. We refer to Stanley \cite{stanley2009} for a thorough exposition of these operators and their basic properties.

For a linear extension $l$ of finite poset $P$ with $|P|=n$, we compute the \textbf{promotion of $l$}, denoted $\delta l$, via the following algorithm.

\begin{enumerate}
    \item We find the element $x\in P$ with $l(x)=1$ and remove its label, leaving $x$ empty.
    \item We then look at the elements covering $x$. We choose the element $y$ with the smallest value of $l$ among them, and slide its value down to $x$ to replace the $1$. This leaves $y$ empty.
    \item Similarly, we slide down to $y$ the smallest value of $l$ among the elements covering $y$, and continue this process, leaving the largest element in the resulting maximal chain without a label. We then label this element $n+1$.
    \item Next, we subtract $1$ from all the labels. The resulting linear extension is $\delta l$, the promotion of $l$.
\end{enumerate}  
 
We then fix the largest value $n$ of $\delta l$. Then we perform promotion on the linear extension of $P$ without the element assigned $n$. We then fix the resulting value of $n-1$. Then we repeat this process until all values have been fixed. The resulting linear extension is called \textbf{the evacuation of $l$}, denoted $\varepsilon l$. A first basic property is that evacuation is an involution on linear extensions. \cref{fig:promotion} gives an example of promotion of a linear extension and \cref{fig:evacuation} continues computing the evacuation of the same linear extension.

 \begin{figure}[H]
\begin{center}
    \begin{tikzpicture}[very thick,scale=0.8]
        \node at (0,0) {};
        \draw [fill] (0,0) circle [radius=0.05];
        \node [below] at (0,0) {$1$};
        \node at (-1,1) {};
        \draw [fill] (-1,1) circle [radius=0.05];
        \node [left] at (-1,1) {$2$};
        \node at (1,1) {};
        \draw [fill] (1,1) circle [radius=0.05];
        \node [right] at (1,1) {$3$};
        \node at (-1,2) {};
        \draw [fill] (-1,2) circle [radius=0.05];
        \node [left] at (-1,2) {$4$};
        \node at (1,2) {};
        \draw [fill] (1,2) circle [radius=0.05];
        \node [right] at (1,2) {$5$};
        \node at (1,3) {};
        \draw [fill] (1,3) circle [radius=0.05];
        \node [right] at (1,3) {$6$};
        \draw (0,0) -- (-1,1) -- (-1,2);
        \draw (0,0) -- (1,1) -- (1,2) -- (1,3);
        \node at (3,1.5) {$\to$};
    \end{tikzpicture}
    \hspace{0.1cm}
    \begin{tikzpicture}[very thick,scale=0.8]
        \node at (0,0) {};
        \draw [fill] (0,0) circle [radius=0.05];
        \node [below] at (0,0) {};
        \node at (-1,1) {};
        \draw [fill] (-1,1) circle [radius=0.05];
        \node [left] at (-1,1) {$2$};
        \node at (1,1) {};
        \draw [fill] (1,1) circle [radius=0.05];
        \node [right] at (1,1) {$3$};
        \node at (-1,2) {};
        \draw [fill] (-1,2) circle [radius=0.05];
        \node [left] at (-1,2) {$4$};
        \node at (1,2) {};
        \draw [fill] (1,2) circle [radius=0.05];
        \node [right] at (1,2) {$5$};
        \node at (1,3) {};
        \draw [fill] (1,3) circle [radius=0.05];
        \node [right] at (1,3) {$6$};
        \draw (0,0) -- (-1,1) -- (-1,2);
        \draw (0,0) -- (1,1) -- (1,2) -- (1,3);
        \node at (3,1.5) {$\to$};
    \end{tikzpicture}
    \hspace{0.1cm}
    \begin{tikzpicture}[very thick,scale=0.8]
        \node at (0,0) {};
        \draw [fill] (0,0) circle [radius=0.05];
        \node [below] at (0,0) {$2$};
        \node at (-1,1) {};
        \draw [fill] (-1,1) circle [radius=0.05];
        \node [left] at (-1,1) {$4$};
        \node at (1,1) {};
        \draw [fill] (1,1) circle [radius=0.05];
        \node [right] at (1,1) {$3$};
        \node at (-1,2) {};
        \draw [fill] (-1,2) circle [radius=0.05];
        \node [left] at (-1,2) {};
        \node at (1,2) {};
        \draw [fill] (1,2) circle [radius=0.05];
        \node [right] at (1,2) {$5$};
        \node at (1,3) {};
        \draw [fill] (1,3) circle [radius=0.05];
        \node [right] at (1,3) {$6$};
        \draw (0,0) -- (-1,1) -- (-1,2);
        \draw (0,0) -- (1,1) -- (1,2) -- (1,3);
        \node at (3,1.5) {$\to$};
    \end{tikzpicture}
    \hspace{0.05cm}
 \begin{tikzpicture}[very thick,scale=0.8]
        \node at (0,0) {};
        \draw [fill] (0,0) circle [radius=0.05];
        \node [below] at (0,0) {$2$};
        \node at (-1,1) {};
        \draw [fill] (-1,1) circle [radius=0.05];
        \node [left] at (-1,1) {$4$};
        \node at (1,1) {};
        \draw [fill] (1,1) circle [radius=0.05];
        \node [right] at (1,1) {$3$};
        \node at (-1,2) {};
        \draw [fill] (-1,2) circle [radius=0.05];
        \node [left] at (-1,2) {$7$};
        \node at (1,2) {};
        \draw [fill] (1,2) circle [radius=0.05];
        \node [right] at (1,2) {$5$};
        \node at (1,3) {};
        \draw [fill] (1,3) circle [radius=0.05];
        \node [right] at (1,3) {$6$};
        \draw (0,0) -- (-1,1) -- (-1,2);
        \draw (0,0) -- (1,1) -- (1,2) -- (1,3);
        \node at (3,1.5) {$\to$};
    \end{tikzpicture}
    \hspace{0.1cm}
 \begin{tikzpicture}[very thick,scale=0.8]
        \node at (0,0) {};
        \draw [fill] (0,0) circle [radius=0.05];
        \node [below] at (0,0) {$2-1$};
        \node at (-1,1) {};
        \draw [fill] (-1,1) circle [radius=0.05];
        \node [left] at (-1,1) {$4-1$};
        \node at (1,1) {};
        \draw [fill] (1,1) circle [radius=0.05];
        \node [right] at (1,1) {$3-1$};
        \node at (-1,2) {};
        \draw [fill] (-1,2) circle [radius=0.05];
        \node [left] at (-1,2) {$7-1$};
        \node at (1,2) {};
        \draw [fill] (1,2) circle [radius=0.05];
        \node [right] at (1,2) {$5-1$};
        \node at (1,3) {};
        \draw [fill] (1,3) circle [radius=0.05];
        \node [right] at (1,3) {$6-1$};
        \draw (0,0) -- (-1,1) -- (-1,2);
        \draw (0,0) -- (1,1) -- (1,2) -- (1,3);
        \node at (3,1.5) {$\to$};
    \end{tikzpicture}
    \hspace{0.1cm}
 \begin{tikzpicture}[very thick,scale=0.8]
        \node at (0,0) {};
        \draw [fill] (0,0) circle [radius=0.05];
        \node [below] at (0,0) {$1$};
        \node at (-1,1) {};
        \draw [fill] (-1,1) circle [radius=0.05];
        \node [left] at (-1,1) {$3$};
        \node at (1,1) {};
        \draw [fill] (1,1) circle [radius=0.05];
        \node [right] at (1,1) {$2$};
        \node at (-1,2) {};
        \draw [fill] (-1,2) circle [radius=0.05];
        \node [circle,draw,left,red] at (-1,2) {$6$};
        \node at (1,2) {};
        \draw [fill] (1,2) circle [radius=0.05];
        \node [right] at (1,2) {$4$};
        \node at (1,3) {};
        \draw [fill] (1,3) circle [radius=0.05];
        \node [right] at (1,3) {$5$};
        \draw (0,0) -- (-1,1) -- (-1,2);
        \draw (0,0) -- (1,1) -- (1,2) -- (1,3);
    \end{tikzpicture}
   \caption{Promotion of the linear extension in \cref{fig;poset stuff examples}.}\label{fig:promotion}
\end{center}
\end{figure}
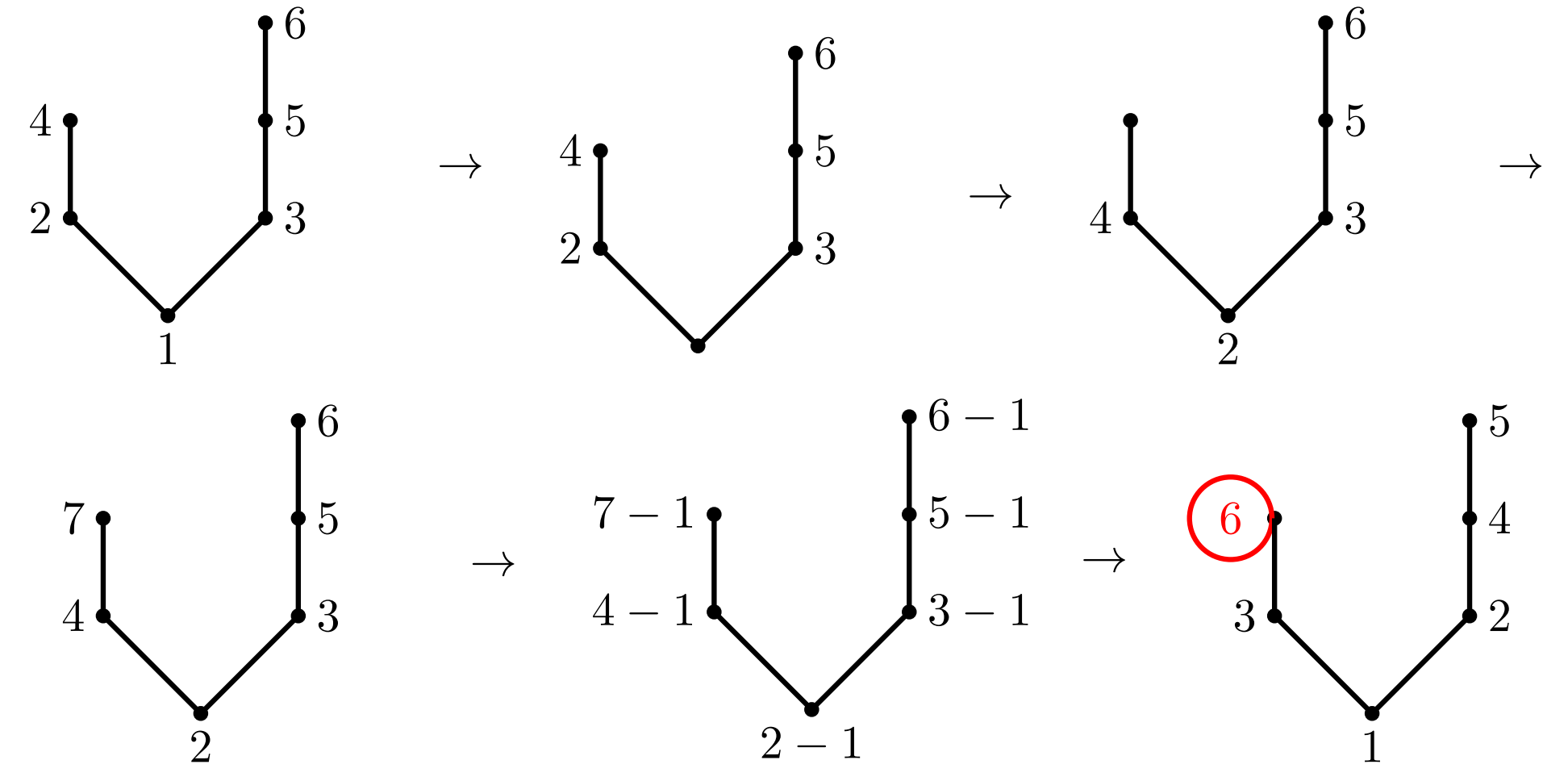
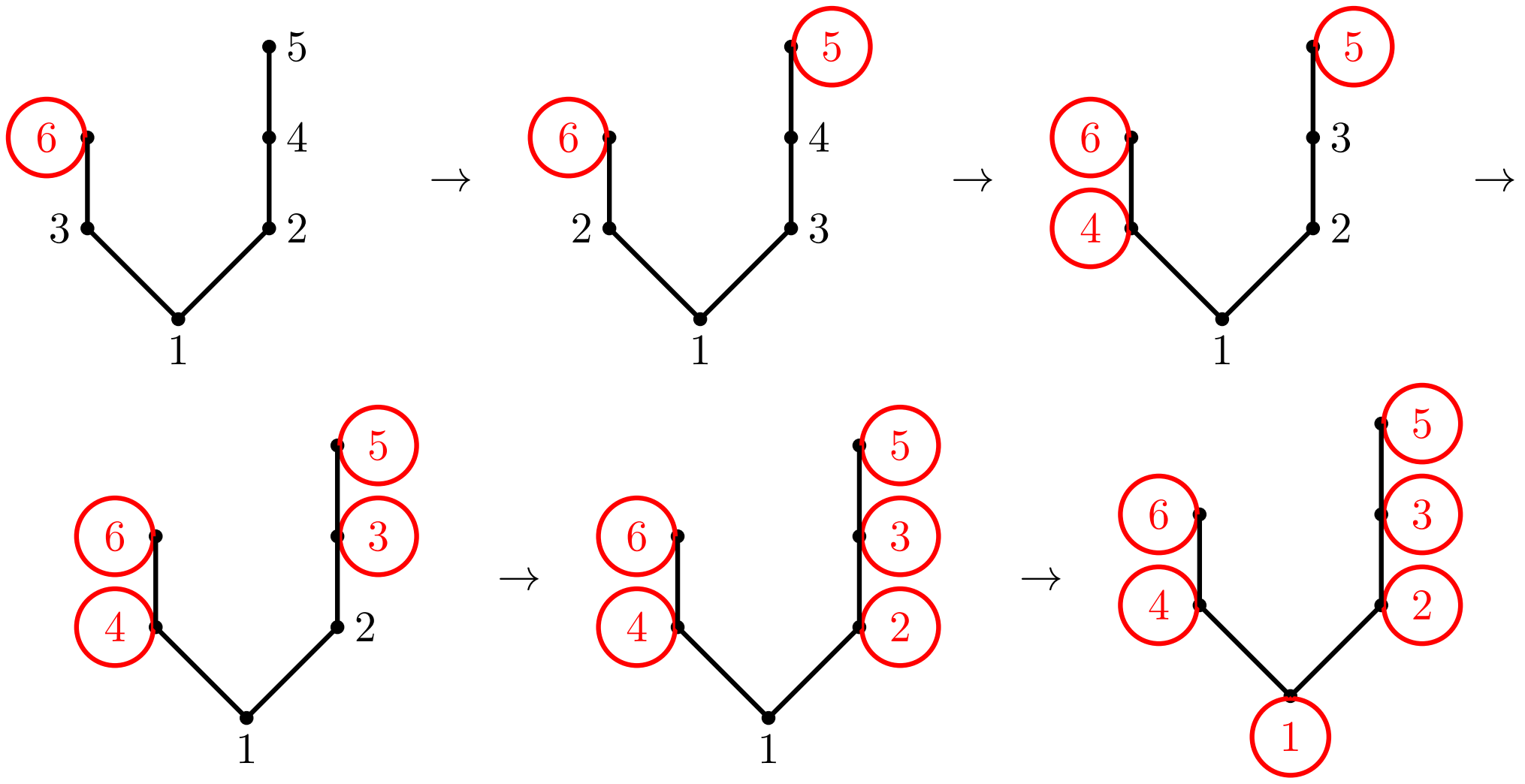
\begin{figure}[H]
\begin{center}
    \begin{tikzpicture}[very thick,scale=0.8]
        \node at (0,0) {};
        \draw [fill] (0,0) circle [radius=0.05];
        \node [below] at (0,0) {$1$};
        \node at (-1,1) {};
        \draw [fill] (-1,1) circle [radius=0.05];
        \node [left] at (-1,1) {$3$};
        \node at (1,1) {};
        \draw [fill] (1,1) circle [radius=0.05];
        \node [right] at (1,1) {$2$};
        \node at (-1,2) {};
        \draw [fill] (-1,2) circle [radius=0.05];
        \node [circle,draw,left,red] at (-1,2) {$6$};
        \node at (1,2) {};
        \draw [fill] (1,2) circle [radius=0.05];
        \node [right] at (1,2) {$4$};
        \node at (1,3) {};
        \draw [fill] (1,3) circle [radius=0.05];
        \node [right] at (1,3) {$5$};
        \draw (0,0) -- (-1,1) -- (-1,2);
        \draw (0,0) -- (1,1) -- (1,2) -- (1,3);
        \node at (3,1.5)
        {$\to$};
    \end{tikzpicture}
    \hspace{0.1cm}
    \begin{tikzpicture}[very thick,scale=0.8]
        \node at (0,0) {};
        \draw [fill] (0,0) circle [radius=0.05];
        \node [below] at (0,0) {$1$};
        \node at (-1,1) {};
        \draw [fill] (-1,1) circle [radius=0.05];
        \node [left] at (-1,1) {$2$};
        \node at (1,1) {};
        \draw [fill] (1,1) circle [radius=0.05];
        \node [right] at (1,1) {$3$};
        \node at (-1,2) {};
        \draw [fill] (-1,2) circle [radius=0.05];
        \node [circle,draw,left,red] at (-1,2) {$6$};
        \node at (1,2) {};
        \draw [fill] (1,2) circle [radius=0.05];
        \node [right] at (1,2) {$4$};
        \node at (1,3) {};
        \draw [fill] (1,3) circle [radius=0.05];
        \node [circle,draw,right,red] at (1,3) {$5$};
        \draw (0,0) -- (-1,1) -- (-1,2);
        \draw (0,0) -- (1,1) -- (1,2) -- (1,3);
        \node at (3,1.5)
        {$\to$};
    \end{tikzpicture}
     \hspace{0.1cm}
    \begin{tikzpicture}[very thick,scale=0.8]
        \node at (0,0) {};
        \draw [fill] (0,0) circle [radius=0.05];
        \node [below] at (0,0) {$1$};
        \node at (-1,1) {};
        \draw [fill] (-1,1) circle [radius=0.05];
        \node [circle,draw,left,red] at (-1,1) {$4$};
        \node at (1,1) {};
        \draw [fill] (1,1) circle [radius=0.05];
        \node [right] at (1,1) {$2$};
        \node at (-1,2) {};
        \draw [fill] (-1,2) circle [radius=0.05];
        \node [circle,draw,left,red] at (-1,2) {$6$};
        \node at (1,2) {};
        \draw [fill] (1,2) circle [radius=0.05];
        \node [right] at (1,2) {$3$};
        \node at (1,3) {};
        \draw [fill] (1,3) circle [radius=0.05];
        \node [circle,draw,right,red] at (1,3) {$5$};
        \draw (0,0) -- (-1,1) -- (-1,2);
        \draw (0,0) -- (1,1) -- (1,2) -- (1,3);
        \node at (3,1.5)
        {$\to$};
    \end{tikzpicture}
     \hspace{0.1cm}
    \begin{tikzpicture}[very thick,scale=0.8]
        \node at (0,0) {};
        \draw [fill] (0,0) circle [radius=0.05];
        \node [below] at (0,0) {$1$};
        \node at (-1,1) {};
        \draw [fill] (-1,1) circle [radius=0.05];
        \node [circle,draw,left,red] at (-1,1) {$4$};
        \node at (1,1) {};
        \draw [fill] (1,1) circle [radius=0.05];
        \node [right] at (1,1) {$2$};
        \node at (-1,2) {};
        \draw [fill] (-1,2) circle [radius=0.05];
        \node [circle,draw,left,red] at (-1,2) {$6$};
        \node at (1,2) {};
        \draw [fill] (1,2) circle [radius=0.05];
        \node [circle,draw,right,red] at (1,2) {$3$};
        \node at (1,3) {};
        \draw [fill] (1,3) circle [radius=0.05];
        \node [circle,draw,right,red] at (1,3) {$5$};
        \draw (0,0) -- (-1,1) -- (-1,2);
        \draw (0,0) -- (1,1) -- (1,2) -- (1,3);
        \node at (3,1.5)
        {$\to$};
    \end{tikzpicture}
         \hspace{0.1cm}
    \begin{tikzpicture}[very thick,scale=0.8]
        \node at (0,0) {};
        \draw [fill] (0,0) circle [radius=0.05];
        \node [below] at (0,0) {$1$};
        \node at (-1,1) {};
        \draw [fill] (-1,1) circle [radius=0.05];
        \node [circle,draw,left,red] at (-1,1) {$4$};
        \node at (1,1) {};
        \draw [fill] (1,1) circle [radius=0.05];
        \node [circle,draw,right,red] at (1,1) {$2$};
        \node at (-1,2) {};
        \draw [fill] (-1,2) circle [circle,draw,radius=0.05];
        \node [circle,draw,left,red] at (-1,2) {$6$};
        \node at (1,2) {};
        \draw [fill] (1,2) circle [radius=0.05];
        \node [circle,draw,right,red] at (1,2) {$3$};
        \node at (1,3) {};
        \draw [fill] (1,3) circle [radius=0.05];
        \node [circle,draw,right,red] at (1,3) {$5$};
        \draw (0,0) -- (-1,1) -- (-1,2);
        \draw (0,0) -- (1,1) -- (1,2) -- (1,3);
        \node at (3,1.5)
        {$\to$};
    \end{tikzpicture}
    \hspace{0.1cm}
    \begin{tikzpicture}[very thick,scale=0.8]
        \node at (0,0) {};
        \draw [fill] (0,0) circle [radius=0.05];
        \node [circle,draw,below,red] at (0,0) {$1$};
        \node at (-1,1) {};
        \draw [fill] (-1,1) circle [radius=0.05];
        \node [circle,draw,left,red] at (-1,1) {$4$};
        \node at (1,1) {};
        \draw [fill] (1,1) circle [radius=0.05];
        \node [circle,draw,right,red] at (1,1) {$2$};
        \node at (-1,2) {};
        \draw [fill] (-1,2) circle [radius=0.05];
        \node [circle,draw,left,red] at (-1,2) {$6$};
        \node at (1,2) {};
        \draw [fill] (1,2) circle [radius=0.05];
        \node [circle,draw,right,red] at (1,2) {$3$};
        \node at (1,3) {};
        \draw [fill] (1,3) circle [radius=0.05];
        \node [circle,draw,right,red] at (1,3) {$5$};
        \draw (0,0) -- (-1,1) -- (-1,2);
        \draw (0,0) -- (1,1) -- (1,2) -- (1,3);
    \end{tikzpicture}
    \caption{Evacuation of the linear extension in \cref{fig;poset stuff examples}.}\label{fig:evacuation}
\end{center}
\end{figure}

\subsection{Weak Order on Permutations}

Let $S_n$ be the symmetric group of permutations of $n$. We think of $S_n$ as a Type A Coxeter group in the standard way with generators the adjacent transpositions $s_i=(i,i+1)$ for $1\leq i\leq n-1$. For more on Coxeter groups, we refer to \cite{bjornerbrenitcxgp}. Recall that a \textbf{reduced expression} of $w\in S_n$ is a shortest possible word in the letters $s_i$ whose product is $w$ and the \textbf{length} of $w$ is number of letters in a reduced expression for $w$.

\begin{definition}[Definition ]\label{def:weak order}
    The (right) \textbf{weak order} on $S_n$ is the partial order on $S_n$ with cover relations $u\lessdot w$ precisely when there is some $s_i$ such that $w=us_i$ and the length of $w$ is one more than the length of $u$. 
\end{definition}

Intuitively, the weak order on $S_n$ has downward cover relations given by swapping adjacent entries in one-line notation from decreasing order to increasing order. We now state two lemmas which we will use to prove some of our results in the following section. The first of these is an automorphism of weak order on $S_n$ which is a special case of an automorphism of weak order on a general Coxeter group.

\begin{lemma}[Prop. 3.1.5 \cite{bjornerbrenitcxgp}]\label{reverselemma}
    $\phi: S_n \to S_n $ given by $\phi(\pi)(i)= n- \pi(n-i+1) +1$ is an poset automorphism of weak order.
\end{lemma}

Second we have a structure theorem for the MCDOs of $J(P)$ from a finite poset $P$ and a linear extension of $P$.

\begin{lemma}[Theorem 6.6 \cite{lacina2022maximalchaindescentorders}]\label{isomcdo}
    Let $P$ be a finite poset with linear extension $l$. Then the MCDO $MP_{\lambda_l}$ is isomorphic to an ideal of weak order on $S_{|P|}$ by mapping each maximal chain to its label sequence.
\end{lemma}

We are now ready to prove our main results.

\section{Isomorphisms of Maximal Chain Descent Orders}\label{sec:results}

In this section, we will give several situations in which two linear extensions of the same finite poset induce isomorphic MCDOs. We begin with several straightforward results culminating in the proof that two linear extensions are related by an automorphism of the poset give rise to isomorphic MCDOs. 
\begin{proposition}
    Let $f: P \to P$ be a poset automorphism and let $l$ be a linear extension of $P$. Then $lf = l\circ f$ is a linear extension of $P$.
\end{proposition}

\begin{proof}
    Suppose $x<y$ in $P$. Since $f$ is an automorphism, $f(x)<f(y) $. Then since $l$ is a linear extension of $P$, we have $lf(x)<lf(y)$.
\end{proof}

\begin{proposition}
Let $I$ be an ideal of $P$ and let $f:P\to P$ be a poset automorphism with inverse $f^{-1}$. Then, $f^{-1}(I)=\sett{f^{-1}(x)}{x\in I}$ is an ideal of $P$.
\end{proposition}

\begin{proof}
Suppose $y\in f^{-1}(I)$ and $z\leq y$ in $P$. Then $y=f^{-1}(x)$ for some $x \in I$. Also, since $f^{-1}$ is a bijection from $P$ to $P$, $z=f^{-1}(w)$ for some $w \in P$. Since $f^{-1}$ is a poset isomorphism and $z\leq y$, we have $w\leq x$. Since $I$ is an ideal of $P$ and $x \in I$, this gives $w \in I$. Thus, $f^{-1}(w) \in f^{-1}(I)$, so $z \in f^{-1}(I)$. Therefore, $f^{-1}(I)$ is an ideal of $P$.
\end{proof}

Now we define the map that will be our isomorphism of MCDOs in this case.

\begin{proposition} \label{F bijection}
Let $F: MP_{\lambda_l} \to MP_{\lambda_{lf}}$ be defined as follows:\\ for $$M:\emptyset \subseteq I_1 \subseteq I_2\subseteq \dots \subseteq P$$ a maximal chain in $J(P)$, $$F(M):\emptyset \subseteq f^{-1}(I_1) \subseteq f^{-1}(I_2) \subseteq \dots \subseteq P.$$ Then, $F$ is a bijection.
\end{proposition}

\begin{proof}
 Since $f^{-1}:P\to P$ is a bijection, $F(M)$ is a maximal chain of $J(P)$. Similarly, $G:MP_{\lambda_{lf}}\to MP_{\lambda_l}$ defined by $G(M):\emptyset \subseteq f(I_1) \subseteq f(I_2) \subseteq \dots \subseteq P$ sends maximal chains of $J(P)$ to maximal chains of $J(P)$. Further, $G$ is the inverse of $F$, so $F$ is a bijection.
\end{proof}

\begin{lemma} \label{label sequence preserved}
    Let $M$ be a maximal chain in $J(P)$ of a finite poset $P$. Let $l$ be a linear extension of $P$. Then, $\lambda_l(M)=\lambda_{lf}(F(M))$ where $F$ is the function defined in \cref{F bijection}.
\end{lemma}
\begin{proof}
    Consider two ideals $I$ and $J$ in $M$ such that $I\lessdot J$. Then $J \backslash I =\set{x}$ for some $x \in P$ and $\lambda_l(I\lessdot J)=l(x)$. Now consider that $f(I)^{-1}\lessdot f^{-1}(J)$ and $f^{-1}(I), f^{-1}(J) \in F(M)$. Then $f^{-1}(J)\backslash f^{-1}(I) =\set{f^{-1}(x)}$ and 
\begin{align*}
    \lambda_{lf}(f^{-1}(I)\lessdot f^{-1}(J))
    &= lf(f^{-1}(x)) \\
    &= l(f(f^{-1}(x))\\
    &= l(x).
\end{align*}
Thus, $\lambda_l(M)=\lambda_{lf}(F(M))$.
\end{proof}

Together the previous lemma and proposition prove our first main result, \cref{thm:poset automorphisms give isom mcdos} from the Introduction.

\begin{theorem} \label{thm:technical poset automorphisms give isom mcdos}
Let $P$ be a finite poset with a linear extension $l$. Let $f:P\to P$ be a poset isomorphism. Let $F: MP_{\lambda_l} \to MP_{\lambda_{lf}}$ be the function defined in \cref{F bijection}. Then $F$ is a poset isomorphism.
\end{theorem} 

\begin{proof}
We know that $F$ is a bijection by \cref{F bijection}. So, it suffices to show that $F$ is order preserving. That is, for any two maximal chains $M$ and $M'$ in $J(P)$,  $M \preceq_{\lambda_l} M'\iff F(M)\preceq_{\lambda_{lf}} F(M')$.
For this, it suffices to show that
$$M \precdot_{\lambda_l} M' \iff F(M) \precdot_{\lambda_{lf}} F(M').$$ To do this, it suffices to show that the label sequences of $M$ and $F(M)$ are the same which we proved in \cref{label sequence preserved} and that $M$ and $M'$ differ in one position if and only if $F(M)$ and $F(M')$ differ in one position.

Suppose $M$ and $M'$ differ in one position $j$, so  $$M: \emptyset \lessdot I_1 \lessdot I_2 \lessdot \dots \lessdot I_{j-1}\lessdot I_j \lessdot i_{j+1} \lessdot \dots \lessdot P$$ and $$M': \emptyset \lessdot I_1 \lessdot I_2 \lessdot \dots \lessdot I_{j-1} \lessdot I'_j \lessdot I_{j+1} \lessdot \dots \lessdot P$$ with $I_j \not= I'_j$. Then $$F(M): \emptyset \lessdot f^{-1}(I_1 )\lessdot f^{-1}(I_2) \lessdot \dots \lessdot f^{-1}(I_{j-1}) \lessdot f^{-1}(I_j) \lessdot f^{-1}(I_{j+1}) \lessdot \dots \lessdot f^{-1}(P)$$ and $$F(M'): \emptyset \lessdot f^{-1}(I_1) \lessdot f^{-1}(I_2) \lessdot \dots\lessdot f^{-1}(I_{j-1}) \lessdot f^{-1}(I'_j) \lessdot f^{-1}(I_{j+1}) \lessdot \dots \lessdot f^{-1}(P).$$ Since $f^{-1}$ is a bijection, $f^{-1}(I_j) \not = f^{-1}(I'_j) $. Thus, $F(M)$ and $F(M')$ differ in one position.
 
On the other hand, suppose that $F(M)$ and $F(M')$ differ in position. Then we apply the same argument using the function $G =F^{-1}$ from the proof of \cref{F bijection} and the fact that $f$ is a bijection to show that $M$ and $M'$ differ in one position.

Thus, $F$ is order preserving, and therefore, $F$ is a poset isomorphism.
\end{proof}

We now turn to our second main result, namely that for a disjoint unions of chains and later ordinal sums of such posets, a linear extension and its evacuation induce isomorphic MCDOs. We begin with a formula for evacuation in the case of a disjoint union of chains, essentially reverse the values and reverse the positions. This formula appears in a more general form in \cite{cactusgroup}. We include a proof here for completeness and simplicity.

\begin{theorem} \label{thm:evacuation formula for disjoint chains}
Let finite poset $P$ be a disjoint union of chains $c_1,c_2,...,c_k$ with $c_i:x_i^1\lessdot x_i^2\lessdot \dots \lessdot x_i^{t_i}$. Let $l$ be a linear extension of $P$. Then the evacuation $\varepsilon l$ is given by $$\varepsilon l(x_j^i)= \sum_{r=1}^{k} t_r-l(x_i^{t_i-j+1})+1.$$
\end{theorem}

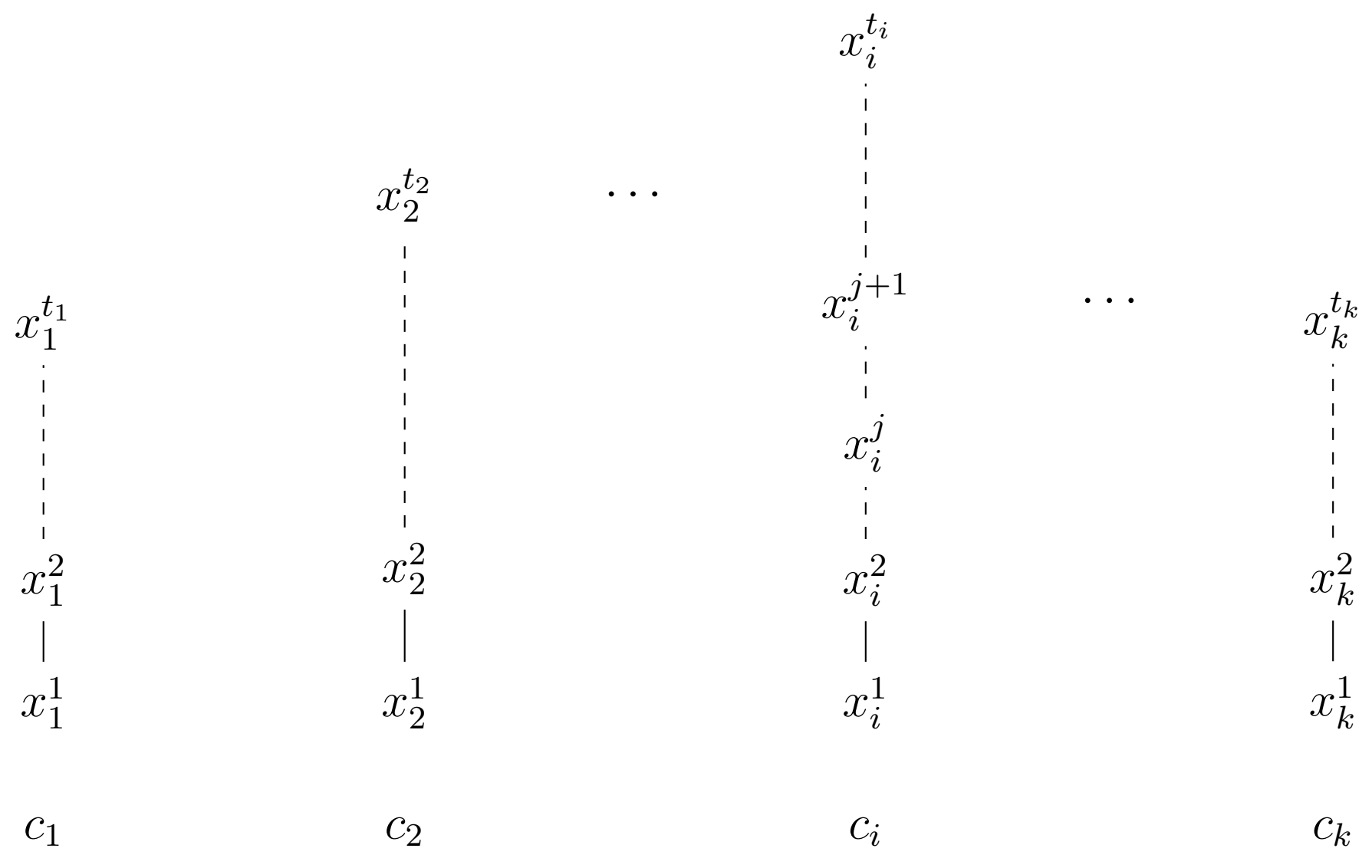
\begin{figure}[H]
\begin{center}
\begin{tikzpicture}[
    node distance=0.35cm,
    every node/.style={font=\large}
]

% Bottom nodes
\node (C1) {$c_1$};
\node (C2) [right=2.5cm of C1] {$c_2$};
\node (Ci) [right=6.5cm of C1] {$c_i$};
\node (Ck) [right=10.5cm of C1] {$c_k$};

% Chain C_1
\node (c1one) [above=0.5cm of C1] {$x_1^1$};
\node (c1two) [above=0.35cm of c1one] {$x_1^2$};
\node (c1top) [above=1.5cm of c1two] {$x_1^{t_1}$};

\draw[] (c1one) -- (c1two);
\draw[dashed] (c1two) -- (c1top);

% Chain C_2
\node (c2one) [above=0.5cm of C2] {$x_2^1$};
\node (c2two) [above=0.45cm of c2one] {$x_2^2$};
\node (c2top) [above=2.5cm of c2two] {$x_2^{t_2}$};

\draw[] (c2one) -- (c2two);
\draw[dashed] (c2two) -- (c2top);

% Dots between C_2 and C_i
\node [right=1.2cm of c2top] {$\cdots$};

% Chain C_i
\node (cione) [above=0.5cm of Ci] {$x_i^1$};
\node (citwo) [above=0.35cm of cione] {$x_i^2$};
\node (cij) [above=0.45cm of citwo] {$x_i^j$};
\node (cijplus) [above=0.45cm of cij] {$x_i^{j+1}$};
\node (citop) [above=1.5cm of cijplus] {$x_i^{t_i}$};

\draw[] (cione) -- (citwo);
\draw[dashed] (citwo) -- (cij);
\draw[dashed] (cij) -- (cijplus);
\draw[dashed] (cijplus) -- (citop);

% Dots between C_i and C_k
\node [right=1.2cm of cijplus] {$\cdots$};

% Chain C_k
\node (ckone) [above=0.5cm of Ck] {$x_k^1$};
\node (cktwo) [above=0.35cm of ckone] {$x_k^2$};
\node (cktop) [above=1.5cm of cktwo] {$x_k^{t_k}$};

\draw[] (ckone) -- (cktwo);
\draw[dashed] (cktwo) -- (cktop);

\end{tikzpicture}
\end{center}
\caption{A poset $P$ that is a disjoint union of chains}
\end{figure}

\begin{proof}
    We will induct on $|P|$. As a base case, let $|P|=1$. We have $P=x_1^1$. Then, $l(x_1^1)=1$ and $\varepsilon l(x_1^1)=1$. Also, $1=1-1+1$.
    
    Suppose the formula holds whenever $|P|=m$. Say $P$ has $|P|=m+1$ with disjoint chains $c_1,c_2,...,c_k$. Let $l$ be any linear extension of $P$. Execute the first promotion step on $l$. Say $l(x_n^1)=1$ for some $n \in \{1,2,...,k\}$. Then, since $P$ is a disjoint union of chains,
    \begin{center}
   
\begin{equation} \label{eq:promotion formula}
 \delta l(x_i^j)= 
\begin{cases}
l(x_i^j)-1, & \text{if } i \neq n,\\
l(x_i^{j+1})-1, & \text{if } i=n \text{ and } 1\leq j < t_n,\\
\displaystyle\sum_{r=1}^{k} t_r, & \text{if } i=n \text{ and } j=t_n.
\end{cases}
\end{equation}
\end{center}
Thus, 
    $$\varepsilon l(x_n^{t_n})=\delta l(x_n^{t_n})=\sum_{r=1}^{k} t_r.$$
    On the other hand, consider
    \begin{align*}
          \sum_{r=1}^{k} t_r-l(x_n^{t_n-t_n+1})+1
          &=\sum_{r=1}^{k} t_r-l(x_n^1)+1\\
          &=\sum_{r=1}^{k} t_r \\
          &=\varepsilon l(x_n^{t_n}).
    \end{align*}
    The remaining promotion steps in the evacuation of $l$ amount to doing evacuation on $P\backslash\{x_n^{t_n}\}$ with linear extension $\delta l$. 

Notice that $P\backslash\{x_n^{t_n}\}$ is also a disjoint union of chains and $\displaystyle |P\backslash\{x_n^{t_n}\}|=\sum_{r=1}^{k} t_r-1=m$ while the chain $c_n$ has one fewer element in $P\backslash\{x_n^{t_n}\}$ than in $P$.  Thus, by the induction hypothesis, if $x_i^{j}\ne x_n^{t_n}$ we have
\begin{center}
\begin{align*}
\varepsilon\delta l(x_i^j)
&=
\begin{cases}
\displaystyle
\sum_{r=1}^{k} t_r-1-\delta l(x_i^{t_i-j+1})+1,
& \text{if } i\neq n,\\[6pt]
\displaystyle
\sum_{r=1}^{k} t_r-1-\delta l(x_i^{t_i-1-j+1})+1,
& \text{if } i=n\text{ and }1\leq j<t_n.
\end{cases}
\\[10pt]
&=
\begin{cases}
\displaystyle
\sum_{r=1}^{k} t_r-( l(x_i^{t_i-j+1})-1),
& \text{if } i\neq n,\\[6pt]
\displaystyle
\sum_{r=1}^{k} t_r-(l(x_i^{t_i-j+1})-1),
& \text{if } i=n\text{ and }1\leq j<t_n.
\end{cases}
\\[10pt]
&=
\begin{cases}
\displaystyle
\sum_{r=1}^{k} t_r- l(x_i^{t_i-j+1})+1,
& \text{if } i\neq n,\\[6pt]
\displaystyle
\sum_{r=1}^{k} t_r-l(x_i^{t_i-j+1})+1,
& \text{if } i=n\text{ and }1\leq j<t_n.
\end{cases}
\end{align*}
\end{center}

We used \cref{eq:promotion formula} at the second equality. Thus, for all $x_i^j$, we have $$\displaystyle \varepsilon l(x_i^j)= \sum_{r=1}^{k} t_r-l(x_i^{t_i-j+1})+1,$$ so the formula holds for all disjoint unions of chains by induction.
\end{proof}

This formula is the key to the proof that evacuation induces an isomorphism of MCDOs in the case of a disjoint union of chains. This is our second main result which appeared in the introduction as \cref{thm:evacuation of disjoint union of chains}.

\begin{theorem}\label{thm:evacuation and disjoint chains}
Let finite poset $P$ be a disjoint union of chains $c_1,c_2,...,c_k$ of the form $c_i:x_i^1\lessdot x_i^2\lessdot \dots \lessdot x_i^{t_i}$.  Then $ MP_{\lambda_l}$ and $MP_{\lambda_{\varepsilon l}}$ are isomorphic.
\end{theorem}
\begin{proof}
    Let $n=t_1+t_2+...+t_k=|P|$. Consider $M: I_0 \lessdot I_1 \lessdot ... \lessdot I_n$.
    For simplicity, represent $M$ as $\mu_1,\mu_2,\dots,\mu_n$ where $\{\mu_i\} = I_i\backslash I_{i-1}$. So, $I_i = \{\mu_1,\mu_2,\dots,\mu_i\} \in J(P)$. Consider a map 
    $$H:MP_{\lambda_l} \to MP_{\lambda_{\varepsilon l}} $$
    defined by $H(M) = \tau_1,\tau_2,...,\tau_n$ where $\tau_i = x_j^{t_j-m+1}$ when $\mu_{n-i+1}=x_j^m$. In spirit, $H$ reverses the chain $M$ and reverses the elements in each chain $c_i$.
    Then $\lambda_l(M) = l(\mu_1),l(\mu_2),\dots, l(\mu_n)$ and $\lambda_{\varepsilon l}(H(M))= \varepsilon l(\tau_1), \varepsilon l (\tau_2),\dots,\varepsilon l (\tau_n)$. Observe that by definition of $\tau_i$ and by \cref{thm:evacuation formula for disjoint chains}
    \begin{align*}
    \varepsilon l (\tau_i)&=\varepsilon l (x^{t_j-m+1}_j) \\
    &= t_1+t_2+...+t_k-l(x_i^{t_j-(t_j-m+1)+1})+1\\
    &= n-l(x_j^m)+1\\
    &= n-l(\mu_{n-i+1})+1.
    \end{align*}
    Thus, $\lambda_{\varepsilon l}(H(M))$ is the result of applying the isomorphism of weak order on $S_n$ given in \cref{reverselemma} to the label sequence $\lambda_l(M)$. Then by \cref{isomcdo}, $H$ is an isomorphism. 
\end{proof}

Our last result has to do with MCDOs derived from the ordinal sum of two posets. We first make several straightforward observations about ordinal sums which we use in the course of the following proofs. We also introduce some notation for later convenience.

\begin{remark}\label{rmk:properties of ordinal sums}
Consider finite posets $P$ and $Q$ with linear extensions $k$ and $l$, respectively. There is a linear extension $k\oplus l$ of the ordinal sum $P\oplus Q$ defined by $\displaystyle k\oplus l(x) =k(x)$ for all $x\in P$ and $\displaystyle k\oplus l(x) =l(x)+|P|$ for all $x\in Q$. Thus, the relative order of values assigned to $P$ by $k\oplus l$ is the same as the relative order of values assigned to $P$ by $k$ and the relative order of values assigned to $Q$ by $k\oplus l$ is the same as the relative order of values assigned to $Q$ by $l$. Next observe that any ideal in $P\oplus Q$ is either an ideal of $P$ or the disjoint union of $P$ with a non-empty ideal of $Q$. Thus, any maximal chain in $J(P\oplus Q)$ is of the form $\emptyset = I_0\lessdot I_1\lessdot \dots \lessdot I_n=P \lessdot J_1\cup P\lessdot J_2\cup P \lessdot \dots \lessdot J_t\cup P=Q\cup P$ where $I_0\lessdot I_1\lessdot \dots \lessdot I_n=P$ is some maximal chain of $J(P)$ and $\emptyset \lessdot J_1\lessdot J_2 \lessdot \dots \lessdot J_t=Q$ is some maximal chain of $J(Q)$. Lastly, for convenience, we will denote a maximal chain in $J(P)$ by $R$ and a maximal chain in $J(Q)$ by $N$. We denote by $R\oplus N$ the maximal chain in  $J(P\oplus Q)$ constructed from $R$ and $N$ as described above.
\end{remark}

\begin{lemma}\label{lem:properties of stacked mcdos}
Let $P$ be a finite poset with a linear extension $k$ and let $Q$ be a finite poset with a linear extension $l$. Let $R$ and $R'$ be a maximal chains of $P$ and $N$ and $N'$ be a maximal chains of $Q$. Then $R\oplus N\precdot_{k\oplus l} R'\oplus N'$ in $M(P\oplus Q)_{k\oplus l}$ if and only if exactly one of the following holds: $R\precdot_{k} R'$ in $MP_k$ and $N=N'$ or $N\precdot_{l} N'$ in $MQ_l$ and $R=R'$.
\end{lemma}

\begin{proof}
    If $R\oplus N\precdot_{k\oplus l} R'\oplus N'$ in $M(P\oplus Q)_{k\oplus l}$, then $R\oplus N$ and $R'\oplus N'$ differ in one spot. Thus, either $R$ and $R'$ differ in one spot and $N=N'$ or $N$ and $N'$ differ in one spot and $R=R'$, but not both. Now consider that $\lambda_{k\oplus l}(R\oplus N)$ lexicographically precedes $\lambda_{k\oplus l}(R'\oplus N')$. In the first case, $\lambda_k(R)$ lexicographically precedes $\lambda_k(R')$ by definition of $k\oplus l$ (see \cref{rmk:properties of ordinal sums}. So, $R\precdot_{k} R'$ in $MP_k$ and $N=N'$. Similarly, in the second case $\lambda_l(N)$ lexicographically precedes $\lambda_l(N')$. So, $N\precdot_{l} N'$ in $MQ_l$ and $R=R'$.

    On the other hand, suppose exactly one of $R\precdot_{k} R'$ in $MP_k$ and $N=N'$ or $N\precdot_{l} N'$ in $MQ_l$ and $R=R'$ holds. In either case, we have that $R\oplus N$ and $R'\oplus N'$ differ in one spot. Again by definition of $k\oplus l$ (see \cref{rmk:properties of ordinal sums}), if $\lambda_k(R)$ lexicographically precedes $\lambda_k(R')$ and $N=N'$, then $\lambda_{k\oplus l}(R\oplus N)$ lexicographically precedes $\lambda_{k\oplus l}(R'\oplus N')$. Similarly, if $R=R'$ and $\lambda_l(N)$ lexicographically precedes  $\lambda_l(N')$, then $\lambda_{k\oplus l}(R\oplus N)$ lexicographically precedes $\lambda_{k\oplus l}(R'\oplus N')$. Thus, in either case, $R\oplus N\precdot_{k\oplus l} R'\oplus N'$ in $M(P\oplus Q)_{k\oplus l}$.
\end{proof}

\begin{theorem} \label{thm:stacking preserves mcdo isom}
Let $P$ be a finite poset with linear extensions $k$ and $k'$ such that $MP_{k}$ is isomorphic to $MP_{k'}$. Similarly, let $Q$ be a finite poset with linear extensions $l$ and $l'$ such that $MQ_{l}$ is isomorphic to $MQ_{l'}$. Let $P \oplus Q$ be the ordinal sum of $P$ and $Q$. Then $M(P\oplus Q)_{k\oplus l}$ is isomorphic to $M(P\oplus Q)_{k'\oplus l'}$.
\end{theorem}

\begin{proof}
    Suppose $|P|=n$ and $|Q|=t$. Now let $f:MP_{k_1}\to MP_{k_2}$ and $g:MQ_{l_1}\to MQ_{l_2}$ be isomorphisms. Define $f\oplus g: MP\oplus Q_{k_1\oplus l_1} \to MP\oplus Q_{k_2\oplus l_2}$ by \[f\oplus g(R\oplus N) = f(R) \oplus g(N).\] 
    
    The function $f\oplus g$ is a bijection since $f^{-1}\oplus g^{-1}$ is its inverse function. Lastly, it suffices to show that for maximal chains $R\oplus N$ and $R'\oplus N'$ of $J(P\oplus Q)$, we have $R\oplus N\precdot_{\lambda_{k_1\oplus l_1}} R'\oplus N'$ if and only if $f\oplus g(R\oplus N) \precdot_{\lambda_{k_2\oplus l_2}} f\oplus g(R'\oplus N')$. 
     
     Suppose $R\oplus N\precdot_{\lambda_{k_1\oplus l_1}} R'\oplus N'$. Then by \cref{lem:properties of stacked mcdos}, either $R\precdot_{k_1} R'$ in $MP_{k_1}$ and $N=N'$ or $N\precdot_{l_1} N'$ in $MQ_{l_1}$ and $R=R'$, but not both. In the first case, $f(R)\precdot_{k_2} f(R')$ in $MP_{k_2}$ since $f$ is a poset isomorphism and $g(N)=g(N')$. Thus, $f\oplus g(R\oplus N) \precdot_{\lambda_{k_2\oplus l_2}} f\oplus g(R'\oplus N')$. Similarly in the second case, $g(N)\precdot_{l_2} g(N')$ in $MQ_l$ since $g$ is a poset isomorphism and $f(R)=f(R')$. So, $f\oplus g(R\oplus N) \precdot_{\lambda_{k_2\oplus l_2}} f\oplus g(R'\oplus N')$.

     On the other hand, $f\oplus g(R\oplus N) \precdot_{\lambda_{k_2\oplus l_2}} f\oplus g(R'\oplus N')$ implies that $R\oplus N\precdot_{\lambda_{k_1\oplus l_1}} R'\oplus N'$ by the same argument using \cref{lem:properties of stacked mcdos} and the fact that $f^{-1}$ and $g^{-1}$ are poset isomorphisms. This concludes the proof.
\end{proof}

\begin{corollary}\label{cor:stacking isomorphism}
    Let $P$ and $Q$ be finite posets with linear extensions $k$ and $l$, respectively. Suppose that for the evacuation of $\varepsilon k$ and the evacuation of $\varepsilon l$, $MP_k$ and $MP_{\varepsilon k}$ are isomorphic and $MQ_l$ and $MQ_{\varepsilon l}$ are isomorphic. Then $M(P\oplus Q)_{k\oplus l}$ and $M(P\oplus Q)_{\varepsilon(k\oplus l)}$ are isomorphic.
\end{corollary}

\begin{proof}
    Since every minimal element of $Q$ covers every maximal element of $P$, the evacuation $\varepsilon(k\oplus l)=\varepsilon k\oplus \varepsilon l$. Then \cref{thm:stacking preserves mcdo isom} gives the result.
\end{proof}

\begin{remark}\label{rmk:stacking Sn evac isomorphism}
Further, suppose the isomorphism of MCDOs induced by evacuation is realized as the isomorphism of weak order on a symmetric group from \cref{reverselemma} as in the proof of \cref{thm:evacuation and disjoint chains}. Then the isomorphism of MCDOs of the ordinal sum from \cref{cor:stacking isomorphism} can also be realized as the isomorphism of weak order on a symmetric group from \cref{reverselemma}.
\end{remark}

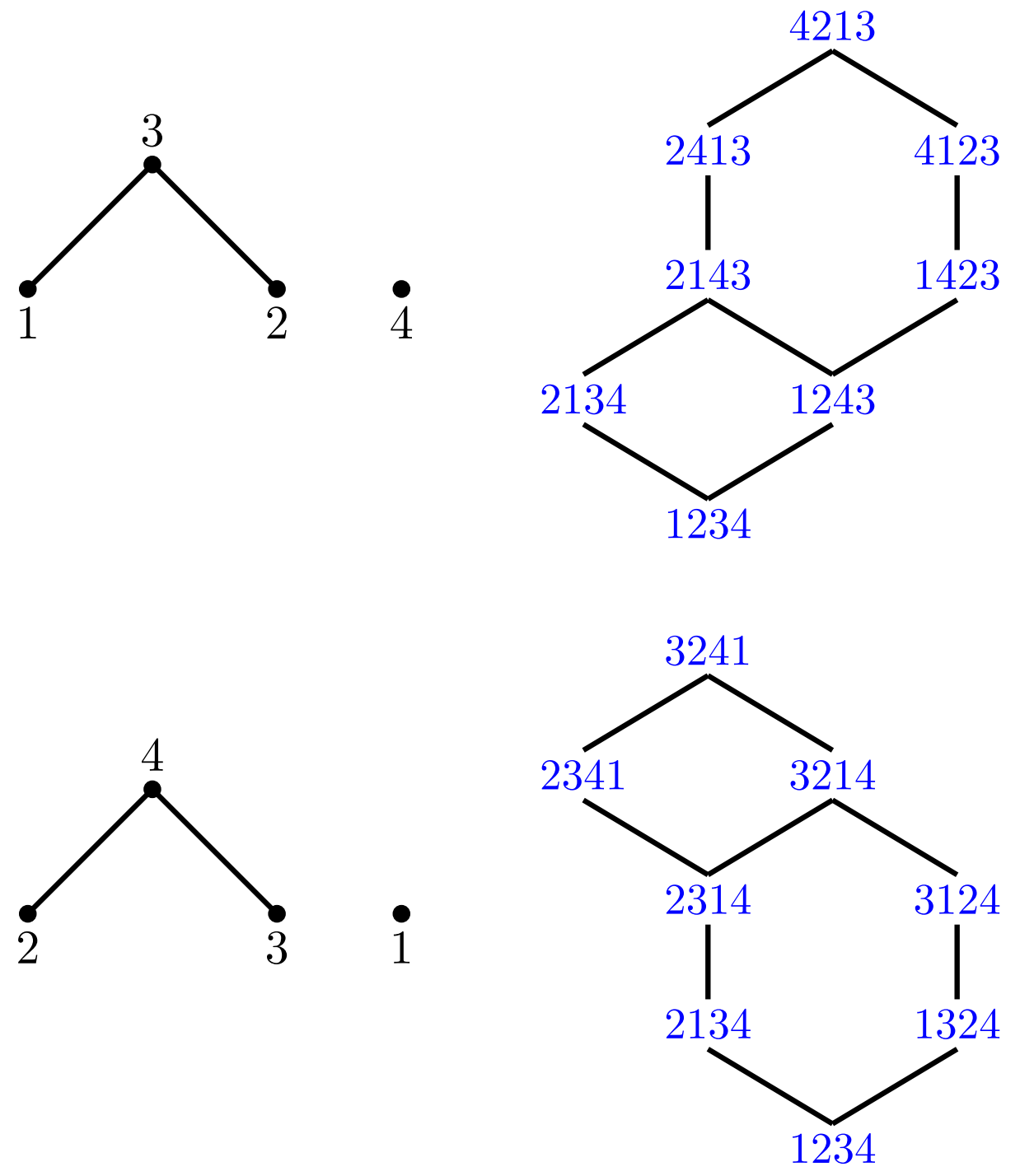
\begin{figure}[h]
    \centering
    \begin{tikzpicture}[very thick,scale=1]
        \node at (-1,0) {};
        \draw [fill] (-1,0) circle [radius=0.05];
        \node [below] at (-1,0) {$1$};
        \node at (0,1) {};
        \draw [fill] (0,1) circle [radius=0.05];
        \node [above] at (0,1) {$3$};
        \node at (1,0) {};
        \draw [fill] (1,0) circle [radius=0.05];
        \node [below] at (1,0) {$2$};
        \node at (2,0) {};
        \draw [fill] (2,0) circle [radius=0.05];
        \node [below] at (2,0) {$4$};
        \draw (-1,0) -- (0,1) -- (1,0);

        \node at (0,-2) {};
    \end{tikzpicture}
    \hspace{0.5cm}
    \begin{tikzpicture}[very thick,scale=1]
        \node[blue] at (0,0) {\small$1234$};
        \node[blue] at (1,1) {\small$1243$};
        \node[blue] at (-1,1) {\small$2134$};
        \node[blue] at (0,2) {\small$2143$};
        \node[blue] at (2,2) {\small$1423$};
        \node[blue] at (0,3) {\small$2413$};
        \node[blue] at (2,3) {\small$4123$};
        \node[blue] at (1,4) {\small$4213$};
    
        \draw (0,0.2) -- (-1,0.8);
        \draw (0,0.2) -- (1,0.8);
        \draw (-1,1.2) -- (0,1.8);
        \draw (1,1.2) -- (0,1.8);
        \draw (1,1.2) -- (2,1.8);
        \draw (0,2.2) -- (0,2.8);
        \draw (0,3.2) -- (1,3.8);
        \draw (2,2.2) -- (2,2.8);
        \draw (2,3.2) -- (1,3.8);
    \end{tikzpicture}
    
    \vspace{0.5cm}
    \begin{tikzpicture}[very thick,scale=1]
        \node at (-1,0) {};
        \draw [fill] (-1,0) circle [radius=0.05];
        \node [below] at (-1,0) {$2$};
        \node at (0,1) {};
        \draw [fill] (0,1) circle [radius=0.05];
        \node [above] at (0,1) {$4$};
        \node at (1,0) {};
        \draw [fill] (1,0) circle [radius=0.05];
        \node [below] at (1,0) {$3$};
        \node at (2,0) {};
        \draw [fill] (2,0) circle [radius=0.05];
        \node [below] at (2,0) {$1$};
        \draw (-1,0) -- (0,1) -- (1,0);

        \node at (0,-2) {};
    \end{tikzpicture}
    \hspace{0.5cm}
    \begin{tikzpicture}[very thick,scale=1]
        \node[blue] at (0,0) {\small$1234$};
        \node[blue] at (1,1) {\small$1324$};
        \node[blue] at (-1,1) {\small$2134$};
        \node[blue] at (-1,2) {\small$2314$};
        \node[blue] at (1,2) {\small$3124$};
        \node[blue] at (0,3) {\small$3214$};
        \node[blue] at (-2,3) {\small$2341$};
        \node[blue] at (-1,4) {\small$3241$};
    
        \draw (0,0.2) -- (-1,0.8);
        \draw (0,0.2) -- (1,0.8);
        \draw (-1,1.2) -- (-1,1.8);
        \draw (1,1.2) -- (1,1.8);
        \draw (1,2.2) -- (0,2.8);
        \draw (-1,2.2) -- (0,2.8);
        \draw (-1,2.2) -- (-2,2.8);
        \draw (0,3.2) -- (-1,3.8);
        \draw (-2,3.2) -- (-1,3.8);
    \end{tikzpicture}
    \caption{A linear extension of a series parallel poset, its evacuation, and the resulting MCDOs.}
    \label{fig:series parallel example}
\end{figure}

\cref{cor:stacking isomorphism} might lead one to wonder if evacuation induces MCDO isomorphisms for series parallel posets, that is, posets which are ordinal sums and disjoint unions of a single element poset. However, this is not the case as illustrated by the simple example in \cref{fig:series parallel example}. This example shows that disjoint union does not necessarily preserve MCDO isomorphisms induced by evacuation.

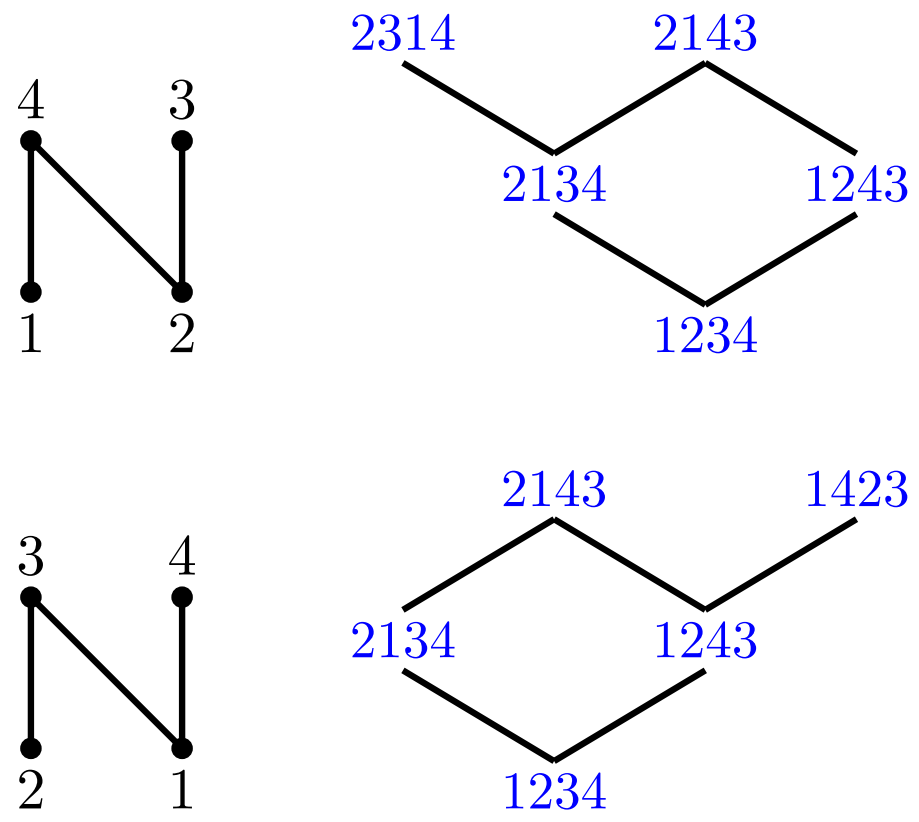
\begin{figure}[H]
    \centering
    \begin{tikzpicture}[very thick,scale=1]
        \node at (-1,0) {};
        \draw [fill] (-1,0) circle [radius=0.05];
        \node [below] at (-1,0) {$1$};
        \node at (-1,1) {};
        \draw [fill] (-1,1) circle [radius=0.05];
        \node [above] at (-1,1) {$4$};
        \node at (0,0) {};
        \draw [fill] (0,0) circle [radius=0.05];
        \node [below] at (0,0) {$2$};
        \node at (0,1) {};
        \draw [fill] (0,1) circle [radius=0.05];
        \node [above] at (0,1) {$3$};
        \draw (-1,0) -- (-1,1) -- (0,0) -- (0,1);
    \end{tikzpicture}
    \hspace{0.5cm}
    \begin{tikzpicture}[very thick,scale=1]
        \node[blue] at (0,0) {\small$1234$};
        \node[blue] at (1,1) {\small$1243$};
        \node[blue] at (-1,1) {\small$2134$};
        \node[blue] at (0,2) {\small$2143$};
        \node[blue] at (-2,2) {\small$2314$};
    
        \draw (0,0.2) -- (-1,0.8);
        \draw (0,0.2) -- (1,0.8);
        \draw (-1,1.2) -- (0,1.8);
        \draw (-1,1.2) -- (-2,1.8);
        \draw (1,1.2) -- (0,1.8);
    \end{tikzpicture}
    
    \vspace{0.5cm}
    \begin{tikzpicture}[very thick,scale=1]
        \node at (-1,0) {};
        \draw [fill] (-1,0) circle [radius=0.05];
        \node [below] at (-1,0) {$2$};
        \node at (-1,1) {};
        \draw [fill] (-1,1) circle [radius=0.05];
        \node [above] at (-1,1) {$3$};
        \node at (0,0) {};
        \draw [fill] (0,0) circle [radius=0.05];
        \node [below] at (0,0) {$1$};
        \node at (0,1) {};
        \draw [fill] (0,1) circle [radius=0.05];
        \node [above] at (0,1) {$4$};
        \draw (-1,0) -- (-1,1) -- (0,0) -- (0,1);
    \end{tikzpicture}
    \hspace{0.5cm}
    \begin{tikzpicture}[very thick,scale=1]
        \node[blue] at (0,0) {\small$1234$};
        \node[blue] at (1,1) {\small$1243$};
        \node[blue] at (-1,1) {\small$2134$};
        \node[blue] at (0,2) {\small$2143$};
        \node[blue] at (2,2) {\small$1423$};
    
        \draw (0,0.2) -- (-1,0.8);
        \draw (0,0.2) -- (1,0.8);
        \draw (-1,1.2) -- (0,1.8);
        \draw (1,1.2) -- (2,1.8);
        \draw (1,1.2) -- (0,1.8);
    \end{tikzpicture}
    \caption{Linear extensions that are neither related by evacuation nor a poset automorphism and the resulting MCDOs.}
    \label{fig:not evac not poset isom}
\end{figure}

Similarly one might wonder if linear extensions that are neither the evacuation of each other nor related by a poset automorphism can induce isomorphic MCDOs. This is the case as illustrated by the two linear extensions and their MCDOs in \cref{fig:not evac not poset isom}. We observe that there are no automorphisms of this poset, and we can simply compute the evacuation of one of these linear extensions and see that it is not the other.

\section{Open Questions}\label{sec:questions}

There are several natural open questions and potential connections that arise from our point of view in this paper. However, we do not include them here to avoid AI further consuming what might be interesting to humans. Please contact the authors directly if you are curious about these questions or other aspects of this work. They are happy to share with fellow mathematicians. 

\section{AI Declaration} AI was not used anywhere in this work, not for mathematical purposes nor writing. All ideas and writing are the authors' own, and they take full responsibility for the correctness of this work.

\bibliographystyle{plainnat}
\bibliography{bibliography}

\end{document}